\documentclass[12pt]{amsart}

\usepackage{amssymb}
\usepackage{amsmath}
\usepackage{amsthm}

\usepackage{amsopn}
\usepackage{algorithm}
\usepackage{algpseudocode}

\usepackage{amsfonts}
\usepackage{graphicx}
\usepackage{nicefrac}
\usepackage{bm}
\usepackage{mathrsfs}
\usepackage{color, colortbl, xcolor}
\usepackage{multirow}
\usepackage[foot]{amsaddr}
\usepackage{fullpage}
\newtheorem{theorem}{Theorem}[section]
\newtheorem{definition}{Definition}[section]

\newtheorem{assumption}[theorem]{Assumption}

\newtheorem{remark}[theorem]{Remark}

\newcommand{\lb}{\left\{}
\newcommand{\rb}{\right\}}
\newcommand{\set}[1]{\left\{ #1 \right\}}

\newcommand{\Def}{\overset{\textbf{def}}{=}}
\newcommand{\RR}{\mathbb{R}}

\newcommand{\PP}{\mathbb{P}}

\newcommand{\rhomax}{\rho_{\text{max}}}

\newcommand{\abs}[1]{\lvert #1 \rvert}
\newcommand{\norm}[1]{\left\| #1 \right\|}

\newcommand{\Hossein}[1]{\textcolor{olive}{#1}}

\newcommand{\tup}[1]{\textup{(}#1\textup{)}}
\usepackage{amsopn}

\newcommand{\BV}{\textbf{BV}}
\DeclareMathOperator{\supp}{Supp}

\newcommand{\pp}{\partial}
\newcommand{\nnorm}[1]{{\left\vert\kern-0.25ex\left\vert\kern-0.25ex\left\vert #1 
    \right\vert\kern-0.25ex\right\vert\kern-0.25ex\right\vert}}

\newcommand{\kernel}{\mathcal K}
\newcommand{\Bor}{\mathscr{B}}

\title{A Nonlocal Degenerate Macroscopic Model of Traffic Dynamics with Saturated Diffusion: Modeling and Calibration Theory}

\author{
    Dawson Do$^{1, *}$ \and
    Hossein Nick Zinat Matin$^{1,*}$ \and
    Masuma Mollika Miti$^{*}$ \and
    Maria Laura Delle Monache $^{*}$
}

\thanks{
    $^{1}$These authors have equal contribution.
}
\thanks{
$^{*}$ Department of Civil and Environmental Engineering, University of California, Berkeley, CA 94720, USA
}
\thanks{
Emails: \{daws, h-matin, masuma\_miti, mldellemonache\}@berkeley.edu
}

\title{A nonlocal degenerate macroscopic model of traffic dynamics with saturated diffusion: modeling and calibration theory}
\begin{document}
\maketitle
\begin{abstract}
In this work, we introduce a novel first-order nonlocal partial differential equation with saturated diffusion to describe the macroscopic behavior of traffic dynamics. We show how the  proposed model is better in comparison with existing models in explaining the underlying driver behavior in real traffic data. In doing so, we introduce a methodology for adjusting the parameters of the proposed PDE with respect to the distribution of real datasets. In particular, we conceptually and analytically elaborate on how such calibration connects the  solution of the PDE to the probability transition kernel proposed by the datasets. 

The performance of the model is thoroughly investigated with respect to several metrics. More precisely, we study the capability of the model in capturing the probability distribution realized by the datasets in the form of the fundamental diagram. We show that the model is capable of approximating the dynamics of the evolution of the probability distribution. To this end, we evaluate the performance of the model with regard to the congestion formation and dissipation scenarios from various datasets. 
\end{abstract}



\textbf{keyword}
traffic flow theory, fundamental diagram, nonlocal, diffusion, partial differential equations






\section{Introduction and related works}
Macroscopic traffic modeling is a key tool that associates a dynamic relationship between traffic flow and speed to the traffic density. The formative \cite{lighthill1955kinematic} and \cite{richards1956shock} kinematic-wave model (known as LWR) reproduces key traffic flow dynamics. Since its introduction, the LWR model has been modified and applied to traffic networks \cite{piccoli2006traffic}, crowd dynamics \cite{rosini2013macroscopic}, and many others; we refer the readers to the survey paper \cite{bellomo2002mathematical}. Variations on the LWR have been developed to study control scenarios such as variable speed limits \cite{holden2015front} and moving bottlenecks \cite{piccoli2009vehicular,delle2014scalar}.

Several reasons contribute to the necessity of improving LWR, suggested first by Lighthill and Whitham themselves \cite{lighthill1955kinematic2}. \cite{nelson2002traveling} lists several discrepancies between empirical traffic measurements and the solution. For example, the LWR model assumes that vehicles change speed instantaneously and based only on local information, neglecting human reaction time and visual perception. Additionally, the LWR model implies there is an injective relationship between density and flow. Plotting observed measurements of measured density and flow on what is known as the fundamental diagram (FD) shows that there is linear pattern up to the critical density, $\rho_c$, corresponding to the maximum flow. However, the distribution of supercritical flows are scattered (Figure \ref{fig:FD_diagram}), implying that multiple values of flow can be associated to higher densities. We refer the interested readers to \cite{nelson2000synchronized, kerner1996experimental} for further observations of the LWR model. As a result of these inconsistencies, research on macroscopic traffic models has evolved in numerous directions.

For example, the macroscopic fundamental diagram (MFD) considers the kinetic model for entire networks, rather than highways, and \cite{geroliminis2008existence} show that the MFD can be reproduced with real-world sensor data. However, the existence of the MFD is only ensured at large-scale and under homogeneous network conditions. For highways, \cite{lighthill1955kinematic2} themselves introduced the first modifications to LWR, which incorporate additional terms in the form of diffusion or viscosity. Diffusive correction is a simple modification to LWR models that preserves the main dynamic properties, while also alleviating the discrepancies with empirical data, such as reproducing scattered FD data. Although initially criticized by \cite{daganzo1995requiem} for implying non-physical dynamics, such as backwards-moving flow, researchers have since produced numerous diffusively-corrected models that ensure physical consistency. The model proposed by \cite{nelson2000synchronized} corresponds to first-order Chapman-Enskog approximation, which includes the continuity equation as well as a diffusively-corrected LWR model and is justified by incorporating reaction time and anticipation distance. Other similar corrections have been proposed in the literature, see for example, \cite{burger2003diffusively, tordeux2018traffic}. Recently, in \cite{matin23Nonlinear}, the authors look at several diffusively-corrected models in the literature and assess both their physical interpretation and their ability to capture empirical results. While diffusive correction modifies the original first-order LWR, second-order models, such as from \cite{payne1971model} and \cite{zhang2000structural}, address the LWR's discrepancies by considering acceleration. The work of \cite{aw2000resurrection} and separately \cite{zhang2002non} addressed criticisms of \cite{daganzo1995requiem}, leading to an extensive literature on the so-called ``ARZ" class of models \cite{seo2017traffic, bellomo2011modeling}. In analogy with fluid dynamics, the first-order models are motivated by mass-balance, while the second-order models incorporate momentum-balance.

Recent works, like the nonlocal LWR, take a separate approach and are motivated by the mathematical interpretation of driver behavior and psychology. Nonlocal terms define velocity on the average density within a ``neighborhood" of each point rather than locally \cite{chiarello2018global, keimer2017existence}. The physical interpretation of this is the drivers' ability to visually assess the traffic ahead of them and react by adjusting their speed \cite{keimer2023nonlocal, blandin2016well, goatin2016well, bayen2022modeling, colombo2012class, huang2022stability}. A review of nonlocal balance law can be found in \cite{keimer2023nonlocal}. While researchers have produced several analytical results regarding the nonlocal LWR, there are few works dedicated to investigating the real-world validity of nonlocal traffic flow models. 

Studying the calibration of the density-flow relationship is essential to traffic flow theory, as empirical data can be used to validate macroscopic models and estimate road capacity, for example. \cite{qu2015fundamental} show that using the weighted least squares method to calibrate freeway flows can accurately fit empirical data in subcritical conditions and mid to high densities. \cite{wang2011logistic} introduces several speed-density functions, while \cite{wang2013stochastic} studies the calibration of stochastic speed-density relationships. Numerous works have been dedicated to the estimation and calibration of the macroscopic fundamental diagram with real data \cite{geroliminis2008existence, dakic2018use, mariotte2020calibration,delle2021three}. Specifically for nonlocal models, in \cite{zhao2024learning} and \cite{huang2024incorporating}, the authors used traffic trajectory data to design physics-informed neural networks to calibrate the fundamental diagram and nonlocal kernel. Their work provides evidence that supports the existence of nonlocal effects. Nonlocal terms have also been explored in cellular automata models \cite{SUN2020132663, SUN2024103083}. The authors show that when look-ahead rules are applied, the data captured from Kinetic Monte Carlo simulations reveal a density-flow relationship commonly observed in real traffic measurements. These recent advances in traffic flow theory indicate that continued investigation into nonlocal models can improve our understanding of macroscopic traffic phenomena.

\subsection*{Contribution}
This paper makes two key contributions to macroscopic traffic flow modeling.  
First, we propose a novel macroscopic traffic flow model that integrates a nonlocal velocity term and a saturated nonlinear diffusive correction. We demonstrate that the proposed model satisfies essential physical principles relevant to traffic flow applications while addressing the limitations of traditional first-order models. Notably, our model accounts for driver perception of traffic dynamics beyond immediate local conditions, offering a more realistic representation of traffic behavior.  

Second, through introducing a novel calibration methodology, we establish that the proposed model effectively captures the probability distribution of the solution implied by empirical density-flow data. This mainly evaluate the performance of the proposed model over a steady-state relationship. In addition, we validate the proposed model's performance in reproducing the transition dynamics of traffic flow across various real-world scenarios. Unlike density-flow distributions, this assessment provides insight into the temporal evolution of the distribution. Based on such capability, we highlight a significant connection between the semigroup solution of the proposed PDE model and the probability transition kernel governing the observed traffic dynamics.  


\section{Mathematical model}
\label{sec:model}
In this paper, we are interested in a nonlocal diffusive partial differential equations with application in traffic flow dynamics, in the form of
\begin{equation}\label{E:main}
    \begin{cases}
        \partial_t \rho(t,x) + \partial_x [\rho(t,x) U([\hat \rho(t, \cdot) * \mathcal K_\gamma](x))] = 0 &, (t, x) \in \mathcal U_T  \\
        \rho(0, x) = \rho_\circ(x) 
    \end{cases}
\end{equation}
where, $\mathcal U_T \Def (0, T) \times \RR$ and for a fixed $\kappa >0$ (the diffusion coefficient), we define the perceived density
\begin{equation} \label{E:fictitious_density}
    \hat \rho \Def \rho + \kappa D(\rho) \Psi(\partial_x \rho),
\end{equation}
with the nonlinear diffusion, 
\begin{equation}\label{E:nonlinear_coeff}
     D(\rho) \Def \rho(1 - \rho), \, \text{for any $\rho \in [0, 1]$}.
\end{equation}
The average velocity function $U(\cdot)$ can be any reasonable decreasing function and $(t, x) \mapsto \rho(t, x) U(t, x)$ represents the flux function. 

The saturation function $u \in \RR \mapsto \Psi(u)$ is considered to be bounded, i.e.,
\begin{equation}\label{E:Psi_bound}
    \norm{\Psi}_{L^\infty(\RR)} \le 1, 
\end{equation}
increasing, odd, and Lipschitz continuous function. 
In particular, 
\begin{equation}\label{E:Psi_function}
    \Psi(s) \to  \Psi_\infty \Def \Psi(\infty), \quad \text{as $s \to \infty$}.
\end{equation}

In practice, such boundedness ensures that (i) the effect of diffusion term is controllable and consequently (ii) the perceived density $\hat \rho$ remains in the valid range $[0, \rhomax]$. Various functions in literature are considered to represent the $\Psi$ function. Depending on the application, some of the functions that have been prevailingly used in the literature are as follows: 
\begin{equation}
    \begin{split}
        \Psi(s) &\Def \tanh(s), \\
        \Psi(s) &\Def \frac{s}{\sqrt{1 + s^2}} \\
        \Psi(s) & \Def \frac{s}{\sqrt{1 + \frac{\nu^2}{c^2} s^2}}
    \end{split}
\end{equation}
where, the first choice of saturate function $\Psi$ is a generic choice which satisfies the requirements. The second definition is a degenerate version of the mean-curvature operator. In the third definition, $\nu>0$ representing the kinematic viscosity and $c>0$ is the maximum speed of propagation (e.g., acoustic). In this paper, we will use $\tanh(\cdot)$. 
\subsection*{Model properties}
We will first look at the perceived density, $\hat \rho$ \eqref{E:fictitious_density}. The nonlocal models, see e.g. \cite{keimer2017existence, keimer2023nonlocal, chiarello2018global} and the references therein, improve the LWR model by considering the weighted average of the density in determining the velocity. In the same spirit, the key intuition for defining the model \eqref{E:main} and \eqref{E:fictitious_density} is the fact that in reality the density felt by the driver might be different from the local density $\rho(t,x)$ and in particular depends on $\partial_x \rho(t,x)$, the changes in the density; see \cite{bonzani2000hydrodynamic, de1999nonlinear}. 
Therefore, we calculate $U(\cdot)$, such that it takes the perceived density into account, i.e. $\hat \rho$. In particular, if $\partial_x \rho$ is positive, drivers feel a larger density as in \eqref{E:fictitious_density}, and similarly for the negative gradient of the density, the drivers feel less density than the density $\rho$. On the other hand, we define the function $\Psi$ which is monotonically increasing, smooth and bounded and hence the partial derivatives are prevented from growing unbounded. In particular, 
\begin{equation}
    \Psi(\partial_x \rho) \to \pm 1 , \quad \text{as $\partial_x \rho \to \pm \infty$}.
\end{equation}
From the traffic flow perspective, such boundedness ensures that the diffusion term is not the dominant term in defining the velocity. Let us next show that from the physical standpoint, $\hat \rho$ is a (perceived) density. 
\begin{theorem}
Let $\hat \rho$ be defined as in \eqref{E:fictitious_density} and $\kappa \in (0, 1)$. Then, $\hat \rho \in [0, 1]$.
\end{theorem}
\begin{proof}
First, it should be noted that when $\partial_x \rho$ increases (decreases) then $\hat \rho$ also increases (decreases) respectively. Equation \eqref{E:Psi_function} implies that as $\partial_x \rho \to - \infty$, the resulting quadratic equation $\rho -\kappa \rho(1 - \rho) - 1$ has the roots $\rho_1 = -\nicefrac{1}{\kappa}$ and $\rho_2 = 1$ and hence for $\rho \in (0, 1)$ and $\kappa \in (0, 1)$ the claim follows. 

For $\partial_x \rho \to \infty$ the quadratic equation $\rho +\kappa \rho(1 - \rho) - 1$ has the roots $\rho_1 = \nicefrac{1}{\kappa}$ and $\rho_2 = 1$ which means for $\kappa \in (0 ,1)$ the claim follows. In fact, the same argument directly can be used for $\abs{\partial_x \rho} < \infty$ by scaling $\kappa \in (0, 1)$ with respect to $\tanh(\partial_x \rho) \in (-1, 1)$. This completes the proof.
\end{proof}
Next, we elaborate on the role of the convolution term. In particular,
\begin{equation}
 (\hat \rho * \mathcal K_\gamma)(t,x) \Def \int_{x}^{x + \gamma} \mathcal K_\gamma(y - x) \hat \rho(t, y) dy  
 \label{e:convol}
\end{equation}
for a kernel $\mathcal K_\gamma$ presents a weighted average of densities $\hat \rho$ in a neighborhood of $x$. The length of this neighborhood is proportional to $\gamma$. In particular, we consider the following assumptions: 
\begin{assumption}
    For a fixed $T >0$, we assume 
    \begin{enumerate}
        \item $\kernel_\gamma$ is decreasing.
        \item $\supp(\kernel_\gamma) \subset [x, x + \gamma]$ \tup{the support of the kernel}. 
        \item We have that $\int_0^\gamma \kernel_\gamma(x) dx = 1$
        \item From practical point of view, in this paper we only consider the previous assumptions on the kernel. From theoretical point of view, we need more restrictions on the kernel to prove the desired properties of the solution. In particular, we need $\mathcal K_\gamma \in  (W^{2, \infty} \cap W^{2,1})(\RR) \cap C^2([0, \gamma])$ \tup{space of weakly differentiable functions} for any fixed $\gamma$. 
        \item The velocity function $U \in C_b^1([0, 1]; \RR_+)$, space of continuously differentiable and bounded functions, which are also decreasing and take values in $[0, U_{\max}]$.      
    \end{enumerate}
\end{assumption}
It can be noted that regardless of the growth rate of $\partial_x \rho$, the flux term 
\begin{equation*}
    \rho U(\hat \rho * \mathcal K_\gamma) \ge 0.
\end{equation*}
This means, by the construction of the problem the flow always moves in the positive direction including the asymptotic case $\partial_x \rho \to \pm \infty$. From the application point of view this implies that the traffic moves in the correct direction. In addition, the model can be expanded to
\begin{equation*}
    \begin{split}
        \rho_t + \partial_x \rho U(\hat \rho * \mathcal K_\gamma) = - \rho \partial_x [\hat \rho * \mathcal K_\gamma] U'(\hat \rho * \mathcal K_\gamma)
    \end{split}
\end{equation*}
where $U'(\cdot)$ denotes the derivative of the function $U$. The term $\partial_x[\hat \rho * \mathcal K_\gamma]$ contains the diffusion term. This implies traffic dynamics are influenced both by the weighted average density and the weighted average derivative of density. From the technical point of view, $U' \le 0$ and hence the diffusion coefficient on the right-hand side will be positive and we have a forward nonlinear parabolic equation. Finally, $D(\rho)$ in the definition of $\hat \rho$ on the right-hand side ensures the diffusion term degenerates properly at $\rho = 0$ and $\rho= 1$. In other words, when $\rho \to 0$ or $\rho \to 1$, the diffusion term vanishes. Consequently, the velocity is not determined by the diffusion in extreme cases and the model does not experience negative velocity. This ensures that the main driving force in the model is advection. 

Before proceeding to the main discussions of this work, we briefly discuss the properties of the solution of the PDE model \eqref{E:main}.
In the case of model \eqref{E:main}, the existence of degeneracy associated with $D(\rho)$ and $\Psi(\partial_x \rho)$ creates discontinuity in the solution. In some cases, the discontinuity will be smoothed out in time. For more theoretical results on the properties of related PDE problems, we refer the interested readers to \cite{campos2021saturated, caselles2013flux, bertsch1992hyperbolic}. 
\begin{definition}[Weak solution]
    A function $\rho \in C((-\infty, T]; L^1 \cap \BV(\RR; \RR))$ is called a weak solution of \eqref{E:main}, if for any test function $\varphi \in C_c^1([0, T] \times \RR; \RR)$, 
\begin{equation} \label{E:weak_general}
\begin{split}
    \int_0^T \int_{\RR}  \rho(t,x) \partial_t \varphi(t, x) dx dt+ \int_0^T \int_{\RR} Q_F(t, x)  \partial_x \varphi(t, x) dx dt - \int_\RR \varphi(0, x) \rho_\circ(x) dx =0
\end{split}\end{equation}
where,
\begin{equation} \label{E:flux} 
Q_F(t, x) \Def \rho(t,x) U(\hat \rho(t) * \kernel_\gamma (x)).
\end{equation}
\end{definition}

\subsection{Numerical scheme}
\label{s:numerical}
In this section, we will briefly describe the numerical scheme used to simulate 
\eqref{E:main}. LDG \cite{cockburn1998local} was chosen for numerically solving the PDE (\ref{E:main}), as it was developed for and is well-suited to the advection-dominated advection-diffusion problem described in section \ref{sec:model}. Considering $Q_F$ as in \eqref{E:flux}, with a slight abuse of notation, we can write (\ref{E:main}) as
\begin{equation}
        \pp_t\rho + \pp_x[Q_F(\rho,\partial_x\rho))] = 0.
\end{equation}
In addition, for any $t \in \RR_+$, we define $R\Def \hat\rho(t,\cdot)*\mathcal{K}_\gamma(\cdot)$, and note that it is a function of $\rho$ and $\pp_x\rho$. Similar to LDG \cite{cockburn1998local}, we split into two equations:
\begin{equation}
    \begin{cases}
        \pp_t\rho + \pp_x[Q_F(\rho,\sigma)] = 0 \\
        \pp_x\rho =\sigma .
    \end{cases}
\end{equation}

We define $\rho_h$ as a function, our approximate solution, within the function space, $V_h$, formed from a set of basis functions, $\{\phi_i^k\}$, for each domain partition $k$. We define $\sigma_h$ similarly, giving:
\begin{align*}
    \rho_h(x)&=\sum_{k=1}^{n}\sum_{i=0}^p \rho_i^k\phi_i^k(x) \\
    \sigma_h(x)&=\sum_{k=1}^{n}\sum_{i=0}^p \sigma_i^k\phi_i^k(x),
\end{align*}
where $\rho_i^k$ and $\sigma_i^k$ are a real-valued coefficients. We use the Legendre polynomials as basis function and chose their coefficients such that they have the property $\phi_i^k(x_i^k) = 1,\,\forall i\in\{0\dots p\}$, where $x_i^k$ are basis points for each domain partition and $p$ is the polynomial degree of our basis. We use Chebyshev nodes as the basis points. Notably, $x_{k-1} = x_p^{k-1}=x_0^k$ and $x_{k} = x_p^{k} = x_0^{k+1}$. 

Under the LDG formulation, our PDE system becomes a set of nonlinear, first-order, coupled ordinary differential equations for each domain partition $k$:
\begin{align}
    &\mathbf{M}^k\mathbf{\sigma}^k = -\mathbf{C}^k\mathbf{\rho}^k+\mathbf{S}^k_1(\rho_h)\label{eq:ldg-1}\\
    &\mathbf{M}^k\dot{\mathbf{\rho}}^k - \mathbf{K}^k(\rho_h,\sigma_h)+\mathbf{S}^k_2(\rho_h)= 0\label{eq:ldg-2}.
\end{align}
From the properties of the selected basis functions, we can view vectors $\dot{\mathbf{\rho}}^k,\mathbf{\rho}^k,\mathbf{\sigma}^k$ as the values of the functions $\dot{\rho}_h,\rho_h,\sigma_h$ at the basis points, $x_i^k$. The other structures are defined as follows: 
\begin{align*}
    M^k_{ij} &= \int_{x_{k-1}}^{x_k} \phi_i^k\phi_j^kdx\\
    C^k_{ij} &= \int_{x_{k-1}}^{x_k} \frac{d \phi_i^k}{dx}\phi_j^kdx\\
    K^k_i(\rho_h,\sigma_h) &= \int_{x_{k-1}}^{x_k}Q_F(\rho_h,\sigma_h)\frac{d \phi_i^k}{dx}dx\\
    \mathbf{S}^k_1(\rho_h) &= \left[ -\rho_0^k, 0, \cdots, 0, \rho_0^{k+1} \right]^\top\\
     \mathbf{S}^k_2(\rho_h) &= \left[-Q(\rho_p^{k-1},\rho_0^{k},R^{k-1}), 0, \cdots, 0, Q(\rho_p^{k},\rho_0^{k+1},R^k)\right]^\top.
\end{align*} 
$Q$ is a numerical flux, in this case, Lax-Friedrich:
\begin{equation*}
    Q(\rho_p^{k},\rho_0^{k+1},R^k) = \frac{1}{2}\left(\left(f\left(\rho_p^{k}\right) + f\left(\rho_0^{k+1}\right)\right)U\left(R^k\right)+\alpha\left(\rho_p^{k}-\rho_0^{k+1}\right)\right),
\end{equation*}
where $R^k=R(x_k)$ and $\alpha = \max\left|\pp_\rho(D(\rho)U(\hat \rho))\right|$. Calculation of $\mathbf{K}^k$, requires the use of Gauss-Legendre quadrature with the number of points, $N_G\geq(p+1)/2$:
\begin{align*}
    K^k_i(\rho_h,\sigma_h) &= \int_{x_{k-1}}^{x_k}Q_F(\rho_h,\sigma_h)\frac{d \phi_i^k}{dx}dx \\
    &= \sum_{g=1}^{N_G}w^k_g\rho_h(x^k_g)U(R(x^k_g))\frac{d \phi_i^k}{dx}(x^k_g)
\end{align*}
where $x^k_g$ and $w^k_g$ are the quadrature points and weights. Calculating $\mathbf{S}^k_2$ and $\mathbf{K}^k$ requires approximating the nonlocal convolution terms for $R^k$ and $R(x_g)$:
\begin{align*}
    R^k &= \int_{x_k}^{x_k+\gamma}\mathcal K_\gamma\left(y-x_k\right)\left(u_h(y)+\kappa D(u_h(y))\Psi(\sigma_h(y))\right)dy\\
    R(x_g) &= \int_{x_g}^{x_g+\gamma}\mathcal K_\gamma\left(y-x_g\right)\left(u_h(y)+\kappa D(u_h(y))\Psi(\sigma_h(y))\right)dy.
\end{align*}
In the case of $R(x_g)$, this needs to be done for each Gauss-Legendre point. The integral over $[x,x+\gamma]$ is computed using piecewise integrals, split at each domain partition. We refer to \cite{chalons2018high} for the detailed formulation of the calculation of the nonlocal convolutions.

At each time step, we solve for $\mathbf{\sigma}^k$ in \eqref{eq:ldg-2} using a linear solver, then it is used to solve for $\dot{\mathbf{\rho}}^k$ in \eqref{eq:ldg-1}. Third Order Runge-Kutta is then used to advance to the next time step. Lastly, a Generalized Slope Limiter and Boundary Preserving limiter is used on the intermediate solutions. Given $\Delta x$, time step $\Delta t$ is chosen such that it satisfies CFL conditions $\frac{\Delta t}{\Delta x} = \beta\frac{1}{(2p+1)\max\left|\pp_\rho Q_F\right|}$ \cite{cockburn2001runge}, where $\beta \leq 1$ is the CFL number. 

 \section{Datasets}

We have used two different datasets to validate the model. The first dataset was obtained from the widely recognized California Department of Transportation Performance Measurement System (PeMS) \cite{choe2002freeway}, known for its comprehensive traffic data. The other data was collected using a drone, providing detailed spatiotemporal traffic data \cite{wu2022decentralized}. 
Each dataset captures different aspects of traffic conditions. The PeMS data primarily comprises of the free-flow regime with occasional congestion, as it is collected through fixed sensors that predominantly measure traffic at highway locations with uninterrupted flow. In contrast, the drone-based dataset captures vehicle trajectories with high spatial and temporal resolution, making it particularly useful for analyzing congestion dynamics, such as stop-and-go waves. By integrating these two datasets, we achieve a more comprehensive calibration process that accounts for both free-flow and congested traffic conditions.
\subsection{PeMS Data}
The data was collected for the entire day from postmiles 20.127 to 24.917 along the I-880 North freeway, covering the regions of Fremont, Union City, and Hayward (CA). The dataset includes information collected through loop detectors installed on various lanes, including mainline, off-ramps, and on-ramps \cite{choe2002freeway}. These detectors recorded average velocity and  traffic flow over a 24-hour period. We used data from April 12, 14, 19, and 21, 2023, with speed and flow measurements recorded at five-minute intervals.  The data was processed to estimate traffic density and its spatial variations over time. 

\subsection{Drone Data}
The data was collected using a single drone over a stretch of highway approximately 800 meters long. The data capture took place during the afternoon rush hours, between 4:00 PM and 6:00 PM, when traffic density was high, particularly under stop-and-go wave scenarios. The drone recorded detailed vehicle movements, capturing the dynamics of traffic flow and interactions.
The data was collected using a DJI Mavic 2 Pro quadcopter, hovering over the roadway to record top-down 4K videos at 30 FPS. This dataset has also been utilized in previous studies to analyze traffic conditions and driver behavior \cite{wu2022decentralized}.


\subsection{Data Preparation}
\label{sub:data_prep}
In this section, we describe how we use the datasets' measurements of observed flow and velocity to estimate the remaining quantities, such as density and its derivative. Let us denote by $\tilde Q(t,x)$ and $\tilde U(t, x)$ the observed flow and velocity at any point $(t,x ) \in \RR_+ \times \RR$, respectively. In addition, $\tilde \rho$ and $\tilde \partial_x \rho$ denote the density and its derivative derived from measured quantities (in the following the discussions, we will refer to them as ``observed"). 
We define
\begin{equation}\label{E:observed_velocity}
    \tilde U(t,x) \Def U(\tilde \rho(t, x),  \tilde \partial_x \rho(t, x))
\end{equation}
Then, considering \eqref{E:main} and \eqref{E:observed_velocity}, we define
\begin{equation}\label{E:observed_flow}
    \tilde Q(t,x) \Def \tilde \rho(t,x) \tilde U(t, x).
\end{equation}
From \eqref{E:observed_velocity} and \eqref{E:observed_flow}, we can retrieve the observed $\tilde \rho(t,x)$ for any $(t,x) \in \RR_+ \times \RR$. In this particular work, we scale the density and speed measurements by their observed maximums, then recompute the flow to normalize the values in the simulation. Additionally, we scale the length of the study section to $1$. It should be noted that scaling is reversible and optional. Let 
\begin{equation}\label{E:observed}
    \mathcal O \Def \set{x_1, \cdots, x_n}
\end{equation}
be the set of observed points. Then, at any time $t$, we use a finite difference approximation to estimate the derivative:
\begin{equation}
     \tilde \partial_x \rho(t, x_j) \Def \frac{\tilde \rho(t,x_{j + 1}) - \tilde \rho(t, x_{j-1})} {x_{j+1} - x_{j-1}}  , \quad j \in \set{2, \cdots, n-1},
\end{equation}
and
\begin{equation}
    \tilde \partial_x \rho(t, x_1) \Def \frac{\tilde \rho(t, x_2)- \tilde \rho(t, x_1)}{x_2 - x_1}, \quad \tilde \partial_x \rho(t, x_n) \Def \frac{\tilde \rho(t, x_n)- \tilde \rho(t, x_{n -1})}{x_n - x_{n - 1}},
\end{equation}
for $j \in \set{1, n}$. We consider $\tilde \partial_x \rho$ as the observed derivative of $\rho$ at any point $(t, x_j)$. We then calculate the empirical perceived density, $\tilde {\hat \rho} $, with \eqref{E:fictitious_density}. Lastly, provided $\mathcal{K}_\gamma$, we approximate the convolution term by constructing a linear interpolation of the empirical perceived density $\tilde{\hat{\rho}}$ profile at every time $t$. Then, we approximate the integral using Guass-Legendre quadrature, only using data points where $x_i < x_n-\gamma$. 

\begin{figure}[h!]
    \centering
    \includegraphics[width=0.45\linewidth]{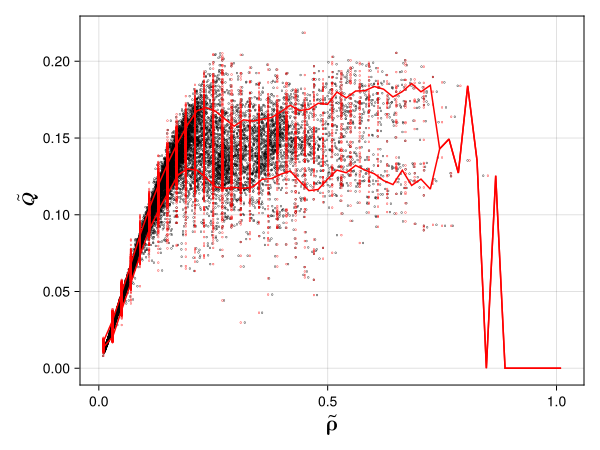}
    \includegraphics[width=0.45\linewidth]{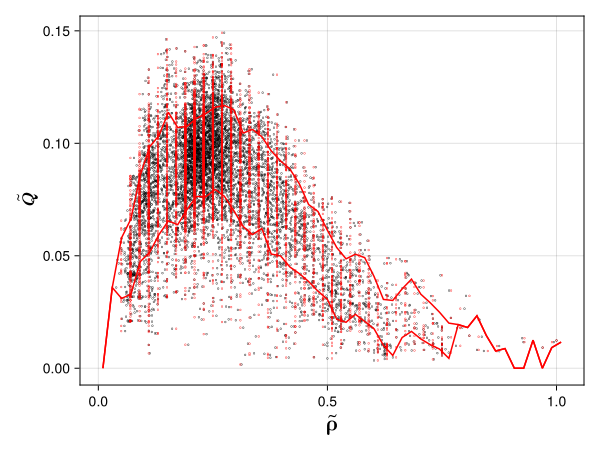}
    \caption{Fundamental diagram of the datasets (left-PeMS, right-drone data). The PeMS dataset mainly contains under-critical density traffic data while the drone data set contain critical and over-critical traffic status. 
    }
    \label{fig:FD_diagram}
\end{figure}

\begin{figure} [h!]
    \centering
    \includegraphics[width=0.45\linewidth]{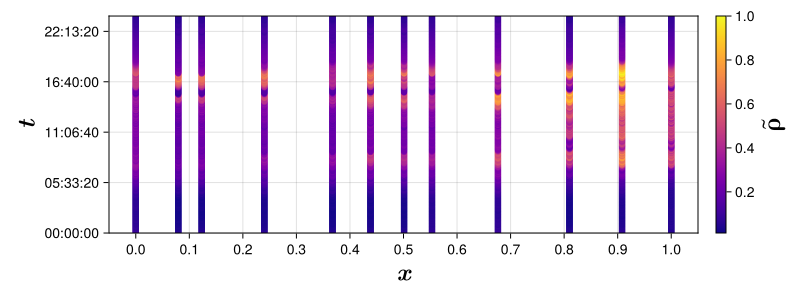}
    \includegraphics[width=0.45\linewidth]{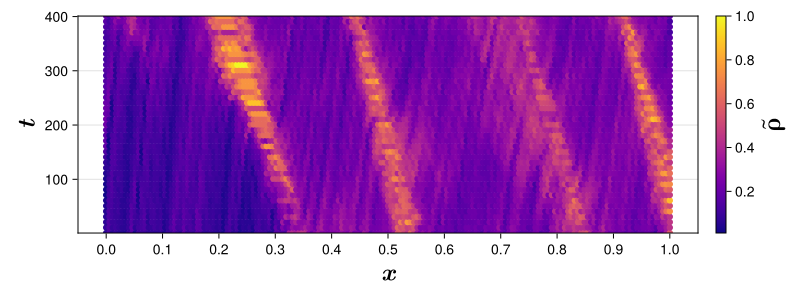}
    
    \caption{Space time plot (Left-PeMS, Right-Drone)}
    \label{fig:x_t_plot_heatmap}
\end{figure}

Figures \ref{fig:FD_diagram} and \ref{fig:x_t_plot_heatmap} depict both datasets in FD and space-time formats. In addition, Figures \ref{fig:drone_3D_plot} and \ref{fig:PeMS_3D_plot} are 3D illustrations of these datasets. 

\section{Methodology}
\subsection{\textbf{Main purpose}} 
Before detailing the calibration methodology, we first establish the conceptual foundation for such calibration by analytically exploring its purpose.
Our goal is to establish a connection between the empirical probability distribution derived from datasets and the PDE solution. To do this, we consider the dataset as a realization of an underlying stochastic process, where the density at different locations follows a probability distribution that evolves over time. The key idea is to describe this evolution using a probability transition semigroup, which captures how the distribution shifts in time. We then show that the semigroup of the PDE solution, which describes the evolution of traffic density in a deterministic setting, is a natural approximation of this stochastic transition. This correspondence justifies calibrating the PDE parameters to align with empirical data.

From a theoretical perspective, we outline how the semigroup framework, when applicable to the solution of the PDE model, can describe the temporal evolution of traffic states as observed in real-world data. 
\begin{remark} [Existence of semigroup solution]
    Proving the existence of such a semigroup for the proposed model \eqref{E:main} falls beyond the scope of this paper. In particular, the effect of the degeneracy $D(\rho)$ and the non-locality on the existence of the semigroup solution of the PDE constitute important topics for future research. Readers interested in similar theoretical developments may refer to \cite{bressan1995semigroup, bressan1996semigroup, andreianov2013semigroup}.
\end{remark}

From a practical perspective, we assess the PDE solution by conducting a statistical comparison with the empirical distribution obtained from the dataset (see Figure \ref{fig:comparison_Q}). Figure \ref{fig:FD_diagram} illustrates the empirical distributions from both the drone and PeMS datasets in the $(\tilde \rho, \tilde Q)$ coordinate space, which serve as the ground-truth probability distribution for our analysis.
\subsubsection*{\textbf{Time evolving probability distribution from the datasets}}
We begin by analyzing the information provided by the dataset. A fundamental assumption is that the empirical data follows an underlying probability distribution that evolves over time. This time-dependent probability distribution, also referred to as the probability transition kernel, describes the stochastic evolution of the traffic state.
(see Definition \ref{def:semigroup} below)

To formalize the concept, we need to define a proper random variable (here the density) over a proper probability space. More accurately, we consider the probability space $(\Omega, (\mathcal F_t)_{t \ge 0}, P)$ equipped with a filtration $(\mathcal F_t)_{t \ge 0}$. From the empirical data, considering the datasets presented in Figure \ref{fig:x_t_plot_heatmap}, at any spatial discretization $\set{x_1, \cdots, x_n}$ as in \eqref{E:observed}, the density $\rho(t,x)$ can be considered as a \textit{random variable} which evolves in time according to the underlying probability distribution of the dataset. More accurately, we consider the random variable
\begin{equation} 
    \bm X(t) = (X_1(t), \cdots, X_n(t)) \Def (\rho(t, x_1), \cdots, \rho(t, x_n)),
\end{equation}
adapted to the filtration $(\mathcal F_t)_{t \ge 0}$ (can be interpreted as information available until time $t$), for $x_i \in \mathcal O$. In this setting, $\bm X(t)$ takes values in the measurable space $(E = \RR^{n}, \Bor(E))$, where $\Bor(E)$ denotes the Borel $\sigma$-algebra (collection of admissible events) on $E$. The density $\rho(t,x)$ is thus analyzed at the location of the collection of observations $\mathcal O$. For abstract definitions of probability space, we refer to classical textbooks such as \cite{durrett2019probability, folland1999real}.
\begin{remark} Let $\Delta x = x_{k +1} - x_k$ for all $k = 1, \cdots, n$, then by letting $\Delta x \to 0$, we can recover the density $\rho(t, x)$ for $x \in \RR$. This transition enables a rigorous analytical discussion, albeit at the cost of introducing additional definitions and notation. \end{remark}
We note that,
\begin{equation}\label{E:random_flow}
\begin{cases}
    \mathfrak d \bm X(t) \Def (\mathfrak d \rho(t, x_1), \cdots, \mathfrak d \rho(t, x_n)), \quad \mathfrak d \rho(t, x_k) \Def  \frac{\rho(t, x_{k + 1}) - \rho(t, x_{k -1})}{x_{k + 1} - x_{k -1}}   \\
   \bm Q(t, x) = \bm Q(\bm X(t), \mathfrak d \bm X(t))
\end{cases}
\end{equation}
i.e., the derivative and the flow functions, respectively, will be the random functions which follow the same probability of $\bm X(t)$. 
In the next section, we discuss the probability transition semigroup associated with the random variable $\bm X(t)$. Some preliminary definitions and concepts can be found in \ref{S:prelim}.
\subsection*{\textbf{Probability measure of the datasets}} 
To facilitate a more structured analysis, we introduce a refined framework for the distribution of the random process $(\bm X(t))_{t \ge0}$. The dataset provides an \textit{empirical probability distribution} of the random process $(\bm X(t))_{t \ge 0}$ which corresponds to a probability transition semigroup $p(t, \xi, A)$ of $(\bm X(t))_{t \ge 0}$, where $ \xi \in \RR^{n}$ and $A \in \Bor(\RR^{n})$ (see \ref{S:prelim} for more details). 
For analytical simplicity of our discussion and without loss of generality, we assume that random process $(\bm X(t))_{t \ge 0}$ follows the \textit{Markov} property with respect to transition semigroup $(p(t,  \cdot, \cdot))_{t \ge 0}$ and the canonical filtration $(\mathcal F_t)_{t \ge0}$. 
\begin{remark}
The discussions in this section can be extended to more general random processes in a straightforward way. In addition, while the Markov assumption may not fully capture all microscopic driver interactions (such as reaction time delays), it serves as a reasonable approximation for macroscopic traffic models where the density at a given time primarily depends on the immediate past rather than a long history.
In other words, the Markov approximation of the random process $(\bm X(t))_{t \ge 0}$ in macroscopic models is reasonable as the proposed PDE model \eqref{E:main} is a first-order model and hence the speed profile is not evolving in time independently (in contrast to the second-order models such as ARZ, \cite{aw2000resurrection}), see \eqref{E:random_flow}. 
\end{remark} 
Mathematically, this means 
\begin{equation}
    P\set{\bm X(t + s) \in A \mid \mathcal F_s} = p(t, \bm X(s), A), 
\end{equation}
which is in fact the definition of a Markov process. 
While this argument is framed within a Markovian context, it can be extended to a broader class of stochastic processes. Notably, by leveraging the Chapman-Kolmogorov identity in the Definition \ref{def:semigroup}, we may conjecture a direct link between the stochastic presentation of traffic states and the deterministic evolution described by a partial differential equation (PDE). Given that the traffic density evolves according to a stochastic transition semigroup, it is natural to ask whether there exists a deterministic PDE semigroup that captures this evolution in an averaged sense. The goal of calibration is to ensure that the semigroup of the PDE solution behaves statistically like the probability transition semigroup observed in real traffic.

Specifically, this connection suggests that there exists a \textit{transition semigroup} $S_t: \mathcal U \to \mathcal U$, $t \in \RR_+$ and with $\mathcal U \subset L^1(\RR)$ that governs the evolution of the traffic density starting from the initial condition $\rho_\circ$. Under this formulation, the solution of PDE \eqref{E:main} can be presented by
\begin{equation}\label{E:semigroup_solution}
    \rho(t, x) = S_t(\rho_\circ)(x). 
\end{equation}
This notation represents the evolution of the traffic density over time using the PDE model.

A rigorous proof of the existence of such a semigroup requires a detailed analysis of the PDE’s properties, which is beyond the scope of this paper. However, its existence plays a crucial role in linking the stochastic formulation of traffic dynamics to the deterministic PDE representation.
\subsubsection*{\textbf{Stochastic-PDE correspondence}}
Since traffic density evolves over time, both in empirical data and in PDE models, a natural question arises: Can the deterministic PDE solution operator $(S_t)_{t \ge0}$ be seen as an approximation of the probability transition semigroup $p(t, \xi, A)$? If so, then the PDE effectively captures the underlying probabilistic traffic dynamics. Our goal is to show that, through appropriate calibration, the PDE solution can statistically approximate the empirical probability evolution.

The preceding discussion suggests that a fundamental approach to demonstrating that the proposed PDE accurately captures the stochastic nature of the empirical dataset is to establish a correspondence between the semigroup $(S_t)_{t \ge 0}$ governing the PDE and the stochastic transition semigroup $(p(t, \cdot, \cdot))_{t \ge 0}$ governing the probability distribution of dataset. To outline such relation, suppose that the semigroup $(S_t)_{t \ge 0}$ exists as in \eqref{E:semigroup_solution}. 

Let $\mathbb F(E)$ denotes the space of right continuous with left limit functions $\alpha: t \in \RR_+ \to E = \RR^n$ equipped with the $\sigma$-algebra $\mathcal F^*$ generated by the maps $\alpha \mapsto \alpha(t)$ (less analytically, we construct the collection of admissible events with respect to the sample space defined by such coordinate maps $\alpha$). Let $\mathbb P$ be the probability measure on $\mathbb F(E)$ which consequently is the law of the Markov process $(\bm X(t))_{t \in \RR_+}$. Let $\gamma(\cdot)$ be the law of initial data $\rho_\circ$. The existence and uniqueness of the probability measure $\PP$ follows directly from the Kolmogorov Extension Theorem and for completeness and accessibility, we provide a brief review of the Kolmogorov Extension Theorem in \ref{S:Kolmo}, which ensures the well-posedness of the probability measure $\PP$ and underpins the validity of the stochastic representation of traffic flow dynamics. This foundational result establishes a rigorous bridge between the empirical stochastic framework and the deterministic PDE model, reinforcing the robustness of the proposed approach in capturing the evolution of traffic states. 
A fundamental result from the theory of Markov processes (see \cite{dynkin1965markov}) states that the finite-dimensional marginals of a stochastic process $(\bm X(t))_{t \in \RR_+}$, denoted by $\PP^U$ on $E^U$, can be explicitly constructed from the transition semigroup $(p(t, \cdot, \cdot))_{t \ge 0}$ and the collection of such finite dimensional marginals satisfies the consistency condition and hence $\PP$ exists. The uniqueness in fact follows from the Chapman-Kolmogorov in Definition \ref{def:semigroup} which describes how probability distributions evolve over time in a Markovian setting. 

This theoretical foundation allows us to conjecture a meaningful relationship between the semigroup of the proposed PDE and the stochastic transition semigroup that governs the evolution of the empirical traffic states. Specifically, if the PDE solution operator
\begin{equation*}
    S_t: \rho_\circ \in L^1(\RR) \to \rho(t) \in \mathbb F(\RR)
\end{equation*}
describes the time evolution of traffic density, we can hypothesize that the probability measure $\PP$ associated with the stochastic process $(\bm X(t)_{t \ge 0})$, which describes the evolution of the traffic density as a stochastic process, can be formally obtained by pushing forward the initial probability distribution $\Phi$ under the PDE semigroup $S_t$. In other words, if $\rho(t,x)$ is the solution to the PDE starting from an initial density $\rho_\circ \sim \Phi$, then the probability distribution of $\rho(t, x)$ at future times is given by a measure
\begin{equation*}
    \PP = S_t \#\Phi.
\end{equation*}
In more intuitive terms, this equation states that the probability distribution of the traffic states at any future time $t$ is obtained by evolving the initial probability distribution $\Phi(\cdot)$ under the deterministic flow induced by the PDE. In other words, for any $B \in \mathcal F^*$, 
\begin{equation*}
    \PP(B) = \Phi(S_t^{-1}(B)).
\end{equation*}
In simple terms, this means that if we evolve the initial density distribution under the PDE model, the resulting density distribution at time $t$ should match the observed probability distribution in traffic data.

As mentioned before, while this paper does not prove the formal existence of the PDE semigroup, prior work on similar traffic flow equations suggests that such a semigroup exists under reasonable conditions. A full rigorous analysis is left for future research.

This formulation suggests a direct connection between the deterministic PDE framework and the stochastic representation of traffic dynamics: if this conjecture holds, then the solution of the proposed PDE \eqref{E:main} naturally corresponds to the Markov process $(\bm X(t))_{t \ge 0}$ inferred from the dataset. 

\begin{figure} [h!]
    \centering
    \includegraphics[width=0.5\linewidth]{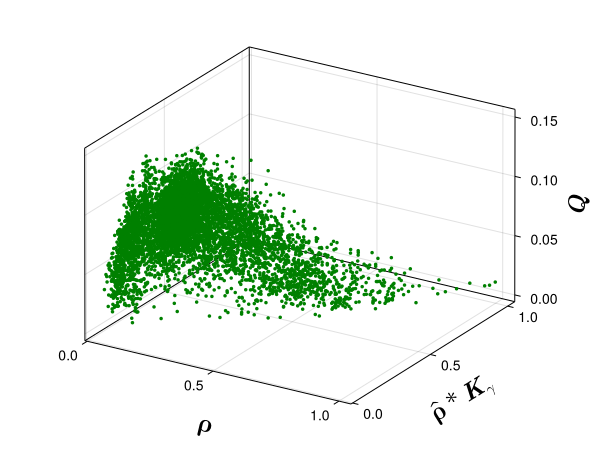}
    \includegraphics[width=0.45\linewidth]{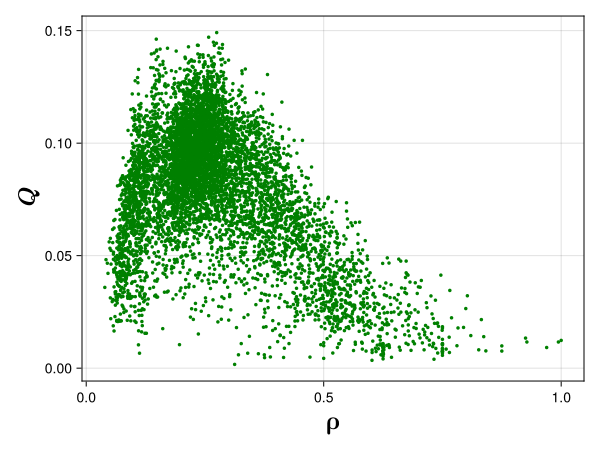}
    \caption{PDE-model prediction of drone data after calibration. Left: 3D plot, and Right: the 2D density-flow plot which matches the real dataset.}
    \label{fig:drone_3D_plot}
\end{figure}

\begin{figure}[ht!]
    \centering
    \includegraphics[width=0.5\linewidth]{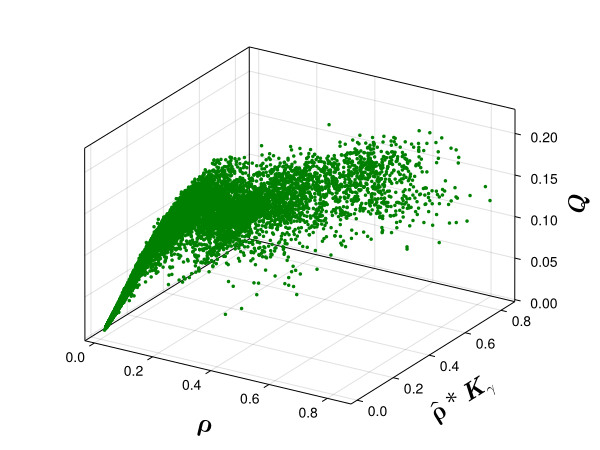}
    \includegraphics[width=0.45\linewidth]{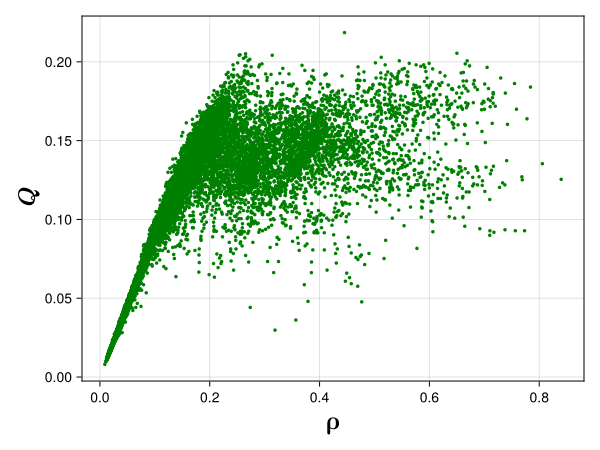}
    \caption{PDE-model prediction of PeMS data after calibration. Left: 3D plot, and Right: the 2D density-flow plot which matches the real dataset.}
    \label{fig:PeMS_3D_plot}
\end{figure}

Figures \ref{fig:drone_3D_plot} and \ref{fig:PeMS_3D_plot} show the prediction of the proposed PDE model \eqref{E:main} after calibration with respect to both datasets.
The left and middle illustrations show 3D plots of $(\rho, \hat \rho * \kernel_\gamma, Q)$ from the PDE model prediction (after calibration). The right illustrations in these Figures show the similarity of the predicted FD to Figure \ref{fig:FD_diagram} of the real datasets. In addition, the discussed correspondence between $(S_t)_{t \ge 0}$ and $(p(t, \cdot, \cdot))_{t \ge 0}$ has been empirically shown and discussed in Figure \ref{fig:comparison_Q}.

\subsection{FD calibration}\label{sub:fd-cal}
As discussed in the previous section, this calibration process is motivated by the treatment of the flux function as a stochastic transition. We start with describing our calibration method in a general framework and then we elaborate on the specific parameters and algorithms used for our problem. Instead of utilizing a least-squares based method, which measures point-wise error, this calibration method measures the similarity between the distributions produced by the empirical data and estimated data. 

Considering \eqref{E:flux}, the Fundamental Diagram (FD) is a function of both $\rho$ and $\partial_x \rho$, unlike the earlier models; see LWR-based model \cite{lighthill1955kinematic, richards1956shock}.
We denote the empirical density and its derivative by $\tilde \rho$ and $ \tilde \partial_x \rho$, respectively. Then, the empirical flow should be approximated by 
\begin{equation}
    \tilde Q(t,x) \approx f(\tilde{\rho}(t,x), \partial_x {\tilde\rho}(t,x);\mathbf{P}),
    \label{eq:k-q-emp}
\end{equation}
where $f(\cdot;\mathbf{P})$ is the approximated density-flow function defined by parameters $\mathbf{P}$. The goal of calibration is to select the appropriate parameters, $\mathbf{P}$ according to the empirical data. In this paper, we will explore several classes of possible functions $f$ and introduce their parameters in the following sections. 

We define the partition 
\begin{equation} 
\set{I_1,I_2,...,I_M} \Def \set{[0,\rho_1),\dots, [\rho_{M-1},\rho_M]}, \quad \rho_M = \rhomax
\end{equation} 
of the interval $[0,\rhomax]$. The data is divided into bins defined by the partition $B_{m}= \set{(\tilde{\rho},\tilde Q)
:\tilde{\rho}\in I_m }$ for all $m\in\{1, \cdots, M\}$. We define $\mu_{\mathbf{P},m},S_{\mathbf{P},m}$ as the average and standard deviation of $\{f(\tilde{\rho}(t,x);\mathbf{P}):\tilde{\rho}\in I_m\}$. 
Similarly, $\mu_{0, m}$ and $S_{0, m}$ represent the mean and standard deviation of the flow values in $\set{(\tilde{\rho},\tilde{Q}) : \tilde{\rho} \in I_m}$, derived from the raw data.
Let $b^{\pm}_{\cdot,m} = \mu_{\cdot,m} \pm S_{\cdot,m}$, the range within one standard deviation of the mean for each bin. We call $\mathbf{b}^{\pm}_{\cdot} = (b^{\pm}_{\cdot,m} \forall m)$, the top and bottom \emph{bands}, respectively, of the fitted data and the empirical measurements. We find $\mathbf{P}$ such that it minimizes the total difference in scatter range:
\begin{equation}
    \mathbf{P}^* = \arg\min_{\mathbf{P}}\left[\sum_{m}\left|b^{+}_{\mathbf{P},m}-b^{+}_{0,m}\right|+\left|b^{-}_{\mathbf{P},m}-b^{-}_{0,m}\right|\right].
    \label{eq:cal-fd}
\end{equation}
For the FD calibration, we also develop two metrics to assess the performance: 1) \textbf{Accuracy}, denoted $\epsilon$, and 2) \textbf{Coverage}, denoted $\Sigma$:
\begin{align}
    \epsilon &= 1 - \frac{\sum_m\left|b^{+}_{\mathbf{P},m}-b^{+}_{0,m}\right|+\left|b^{-}_{\mathbf{P},m}-b^{-}_{0,m}\right|}{\sum_m\max\left(b^{+}_{\mathbf{P},m},b^{+}_{0,m}\right)-\min\left(b^{-}_{\mathbf{P},m},b^{-}_{0,m}\right)} \\
    \Sigma &= \frac{\sum_m\min\left(b^{+}_{\mathbf{P},m},b^{+}_{0,m}\right)-\max\left(b^{-}_{\mathbf{P},m},b^{-}_{0,m}\right)}{\sum_mb^{+}_{0,m}-b^{-}_{0,m}}.
\end{align}

Accuracy, $\epsilon$, describes how similar the bands of the fitted data are to the bands of the measured data. Coverage, $\Sigma$, describes the capability of the fitted data's bands to fully cover the area between the bands of the measured data. 

\subsubsection{Parameter optimization}
For each model, we assume a given $\mathcal K_\gamma$, and search over possible velocity functions and values of $\kappa$ as described in the following sections. We consider two classes of velocity functions, the exponential and piecewise-cubic velocity functions.

\subsubsection*{{\textbf{Exponential velocity function}}} In this paper, we explore the density-speed function explored by \cite{newell1961theory}: $
    U(\rho) \Def v \left(1 - \exp \lb c/v \left(1 - \rhomax/\rho \right) \rb \right)$. There are only two parameters in this speed function, plus $\kappa$, so the optimization can be solved with brute-force.

\subsubsection*{{\textbf{Cubic approximation}}} \label{S:polynomial_aprroximation}
\begin{figure}[H]
    \centering
    \includegraphics[width=0.45\linewidth]{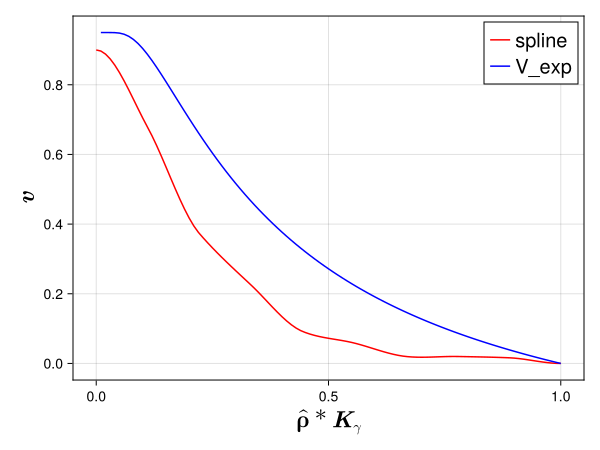}
    \caption{Comparison of Velocity Profile: Exponential and Polynomial}
    \label{fig:velocity}
\end{figure}

In addition, we also search over all decreasing functions represented by cubic spline interpolation, defined by the set control points $p_i$. In this case, the parameters, $\mathbf{P}$, are the position of the control points with which we construct the interpolation. Let $n$ be the number of control points, $p_i = (\rho,v)_i$, $i = 0, \cdots, n-1$, and assume that we fix each control point's density to be equidistant from each other, i.e. $\rho_i = (i-1)\times\rhomax/(n-1)$, $i = 0, \cdots, n-1$. For a suitable search range over $N$ values of speed, using brute-force we would need to search $N^n$ possible cases to find the optimal function value. Figure \ref{fig:velocity} shows both exponential and polynomial approximation of the velocity function. In addition, we noticed that for sufficient fit, $n$ should not be less than 6, which necessitates the iterative optimization Algorithm \ref{alg:spline}. 
We fix $\kappa$ to find the optimal control points for each $\kappa$. First, we evaluate $\hat \rho * \mathcal K_\gamma$ for all $\tilde{\rho}$, giving pairs $(\tilde{\rho}, \hat \rho * \mathcal K_\gamma)$. To select control point $p_i$, we construct the cubic spline interpolation using points up to $p_i$ while including $p_n$, (i.e. $p_0,\dots, p_i,p_n$), for all possible values of $p_i$. The value of $p_i$ is restricted such that the cubic spline interpolation is still decreasing. Then, given this truncated $f$, which is only defined from $\rho\in[\rho_i,\rho_{i+1}]$, we evaluate the statistics $b^{\pm}_{\mathbf{P},m}$, defined above, up to $m$ such that $\rho_m<\rho_i+1$. We then select $p^*_i$ such that it minimizes \eqref{eq:cal-fd}. We repeat this for all $i$. To initialize, we set $v_0$, the speed of the first control point, to be the average speed at $\rho=0$ and $v_n=0$ (Refer to Algorithm \ref{alg:spline} for details). Finally, across all the evaluated $\kappa$, we select $\kappa$ and the corresponding cubic spline interpolation that yielded the overall minimal value for \eqref{eq:cal-fd}. This algorithm searches over the space of $N\times n$, which is a vast improvement over the brute-force algorithm and yields sufficient results. 
\begin{algorithm} 
\caption{Spline Optimization Heuristic}
\label{alg:spline}
\begin{algorithmic}[1]
\State With defined $\kappa, \gamma, \mathcal K_\gamma$, calculate $\hat \rho * \mathcal{K}_\gamma$ as in Section \ref{sub:data_prep}.
\State Set $n$ as the number of control points
\State Initialize $\rho_i = (i-1)\times\rhomax/(n-1)$ for all $i=1\dots n$.
\State Initialize $v_i=0$ for all $i=1\dots n$.
\State Set $v_1 = v_{\max}$, where $v_{\max}$ is the average observed speed of measurements in the first bin, $B_1$ (see Section \ref{sub:fd-cal}).
\State Define $\mathbf{v}_{\text{ops}} = \set{\frac{j}{N}\times v_{\max}|\, j=0\dots N}$ as a range from $0$ to $v_{\max}$ with $N+1$ steps. This is set of options for each $v_i$. 
    
\For{$i = 1$ to $n-1$}
    \State Initialize $d_{\min} = \infty$ to store the minimum difference.
    \State Set $v_\text{cand} = v_i$
    
    \For{each $v\in\mathbf{v}_{\text{ops}}\cap[\frac{n - i}{N}v_{\max},v_i)$}
        \State Set $v_{i+1} = v$
            \State Define $\mathbf{m}_i=\set{ m : I_m \in \set{I_1,I_2,...,I_M}\cap[\rho_i,\rho_{i+1})}$.
            \State Linear interpolate between $(\rho,v)_i$ and $(\rho,v)_{i+1}$ to approximate $U$
            \State Calculate $\{f(\tilde{\rho}(t,x);\mathbf{P}):\tilde{\rho}\in I_m\}\,\forall m\in \mathbf{m}_i$ 
            \State Calculate $b^{\pm}_{\mathbf{P},m}\, \forall m\in \mathbf{m}_i$
            \State Calculate $d = \sum_{m\in\mathbf{m}_i}\left|b^{+}_{\mathbf{P},m}-b^{+}_{0,m}\right|+\left|b^{-}_{\mathbf{P},m}-b^{-}_{0,m}\right|$.
            \If{$d < d_{\min}$}
                \State Update $d_{\min}=d$.
                \State Update $v_\text{cand}=v$.
            \EndIf
    \EndFor
    \State Set the optimal control point value: $v_{i+1}=v_\text{cand}$
\EndFor

\State Set $v_n=0$ enforce the final constraint.
\State Return the cubic spline defined by the set $\set{(\rho,v)_i,\,i=1\dots n}$ (with boundary derivatives enforced to 0.)
\end{algorithmic}
\end{algorithm}


\subsection{Solution-based calibration}\label{sec:solution-based}
In addition to the fundamental diagram calibration, we explore calibration based on the PDE solution. In order to test the ideal performance of the various models, we  used the numerical solver (described in section \ref{s:numerical}) to calibrate $\mathbf{P}$. We select $\tilde{\rho}(t,x)$ for a specific time frame, $t \in [T_0,T_f]$. We set the initial condition to be $\tilde{\rho}(T_0,x)$ and the boundary conditions to be $\tilde{\rho}(t,x_L)$ and $\tilde{\rho}(t,x_R)$ at the ghost boundary cells. We then run the solver until $T_f$, obtaining our solution $\rho$. We find $\mathbf{P}$ such that it minimizes the $L^2$-norm:
\begin{equation}
    \mathbf{P}^* = \arg\min_{\mathbf{P}}\left[\sqrt{\int_{T_0}^{T_f}\int_{x_L}^{x_R}\left(\tilde{\rho}-\rho\right)^2dxdt}\right].
    \label{eq:cal-sol}
\end{equation}

\subsubsection{Parameter optimization}

Similar to before, we assume a given $\mathcal K_\gamma$, and test over the possible parameters for the velocity function and $\kappa$. For the solution-based calibration, we only present the optimization with the Newell density-speed function. Again, the Newell density-speed function can be searched by brute-force to find $\mathbf{P}$. We do not fit a spline velocity function using this calibration because this would be very computationally expensive, given the number of parameters required to find the spline control points. 

\section{Results for Fundamental Diagram-based calibration}
In this section, we explain the calibration process used to fit the model parameters to the empirical FD measurements from each dataset. We have calibrated the parameters of the FD using two distinct datasets: drone-based trajectory data and PeMS sensor data. The trajectory data captures a congestion period with high-resolution, but lacks free-flow information, while the sensor data predominantly samples the free-flow regime. In this work, we use both datasets to show the effectiveness of the calibration across all traffic conditions. 

\begin{figure}[ht!]
    \centering
    \includegraphics[width=0.45\linewidth]{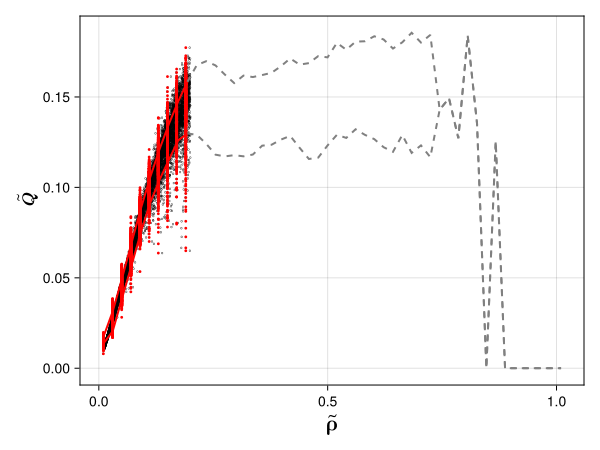}
    \includegraphics[width=0.45\linewidth]{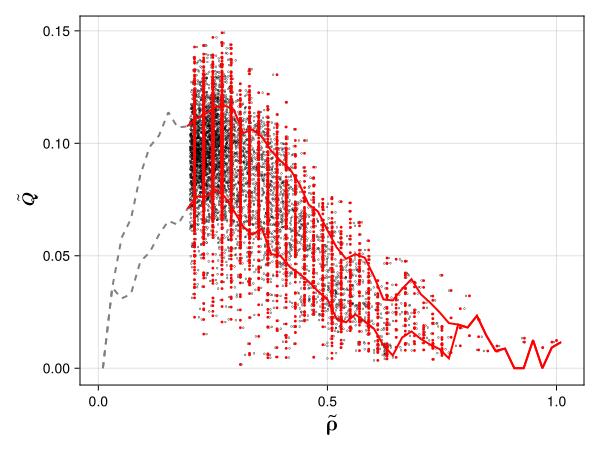}
    \caption{Fundamental diagram (left-free flow (PeMS), right- congested (drone))}
    \label{fig:FD_diagram_PeMS}
\end{figure}

As stated above, we use the drone-collected data to calibrate the FD in congestion. Therefore, we only use samples where $\tilde \rho \geq 0.2 \approx \rho_c$; the critical density. In contrast, we use samples where $\tilde \rho < 0.2$ in the PeMs data, (Figure \ref{fig:FD_diagram_PeMS}). 

\begin{figure}[ht!]
    \centering
    \includegraphics[width=0.45\linewidth]{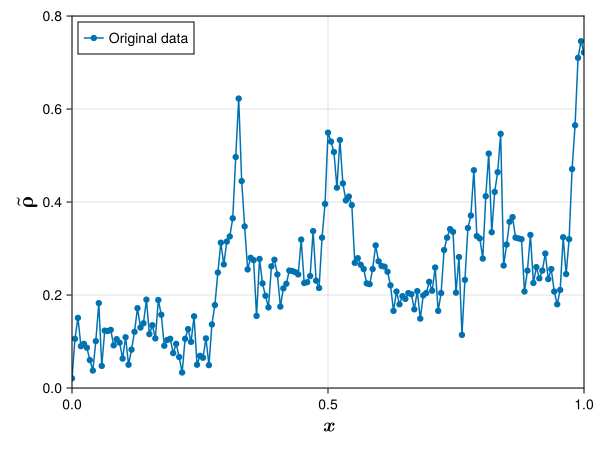}
    \includegraphics[width=0.45\linewidth]{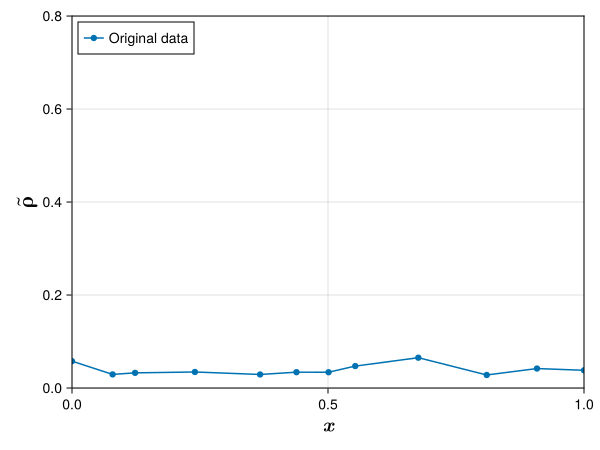}
    \caption{Density plot for the original datasets (left-drone data, right-PeMS data)}
    \label{fig:density_plot}
\end{figure}

\begin{figure}[ht!]
    \centering
    \includegraphics[width=0.5\linewidth]{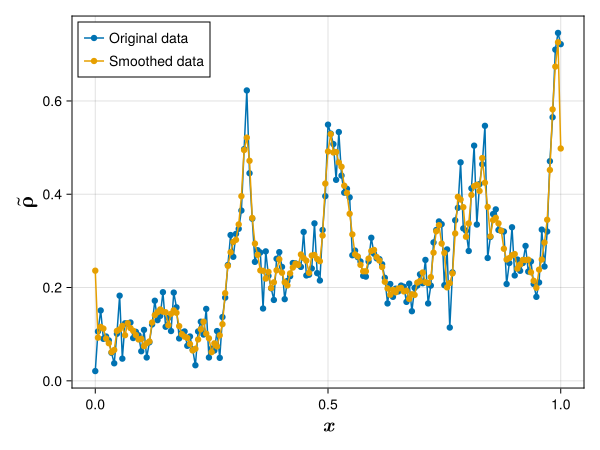}
    \caption{Density profile for the drone data after adjusting density}
    \label{fig:box_filter}
\end{figure}

\begin{figure}[ht!]
    \centering
    \includegraphics[width=0.45\linewidth]{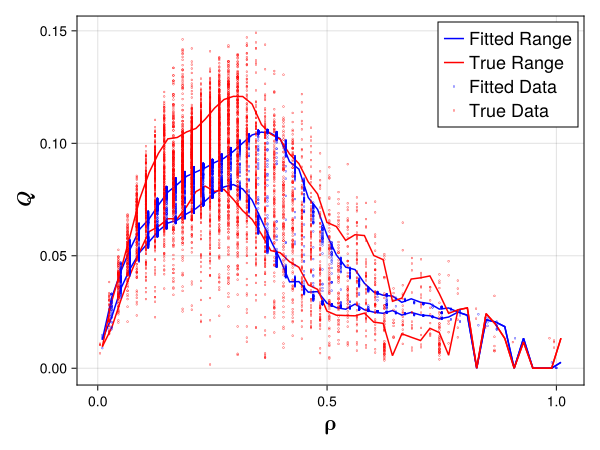}
    \includegraphics[width=0.45\linewidth]{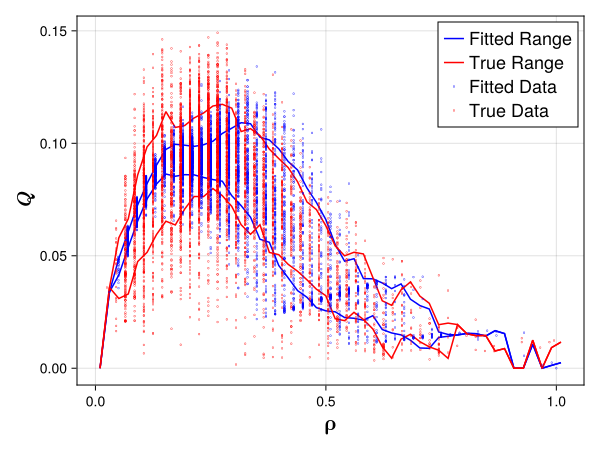}
    \caption{ Fundamental diagram(left-without adjusting density, right-after adjusting density)}
    \label{fig:FD_diagram_K1_K2_K3}
\end{figure}

The first dataset consists of drone-based trajectory data, which is critical for calculating the derivative of density, $\tilde \partial_x \rho$. However, its discretization affects the accuracy of the saturation function $\Psi$ \eqref{E:fictitious_density}. As is, the magnitudes of the estimated empirical $\tilde \partial_x{\rho}$ in the drone-collected dataset can be large (Figure \ref{fig:density_plot}), causing the result of $\Psi(\tilde \partial_x{\rho})$ to reach asymptotic saturation values close to $1$ and $-1$ often. As a result, most values for the estimated perceived density end up polarized near $\rho\pm\kappa D(\rho)$. This produces a bimodal distribution for the approximate flow (i.e.  $f(\tilde{\rho},\tilde \partial_x \rho)$) for a given density, which is dissimilar to the distribution of the empirical flow (Figure \ref{fig:FD_diagram_K1_K2_K3}). Therefore, as an initial pass, we apply a linear box filter of radius 1 to the density profile at every time step (i.e. a moving average of each point with its immediate neighbors) (Figure \ref{fig:box_filter}). Additionally, we also introduce additional parameters ($K_1$, $K_2$, $K_3$) into $\Psi$ to allow the calibration to reduce the polarization effect:
\begin{equation}
\Psi(u) = \tanh\frac{K_1 u - K_2}{K_3}.
\label{eq:phi}
\end{equation}
Applying a linear transformation to $\partial_x \rho$ allows the model to control the magnitudes that will be clamped to $-1$ and $1$.

The PeMS dataset, which primarily captures free-flow traffic conditions, is calibrated without these parameters, as the empirical $\tilde \partial_x \rho$ has relatively small magnitudes. The plot for the free-flow region confirms that PeMS data calibration remains unaffected without these parameters (Figure \ref{fig:calibrated_FD}). 

Figure \ref{fig:density_plot} illustrates the impact of discretization on density changes achieved using $K_1$, $K_2$, and $K_3$. Without these parameters, as shown in the left panel of Figure~\ref{fig:FD_diagram_K1_K2_K3}, the fundamental diagram shows that most of the data points are concentrated along the boundary of the fitted range (i.e. most measurements of $\tilde \partial_x {\rho}$ are clamped).
The calibration with $K_1$, $K_2$, and $K_3$ (right panel of Figure~\ref{fig:FD_diagram_K1_K2_K3}) results in a more natural distribution of flow values. These adjustments allow the calibration process to better capture the distribution of traffic flow, particularly in congestion.

\begin{figure}[ht!]
    \centering
    \includegraphics[width=0.4\linewidth]{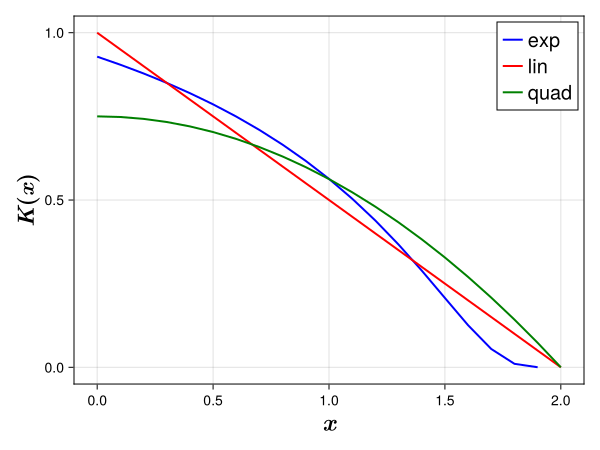}
    \caption{Visualization of Linear, Quadratic, and Exponential Kernels for $\gamma=2$}
    \label{fig:kernel}
\end{figure}

In these experiments, we test three types of kernels for $\mathcal K_\gamma$ \eqref{e:convol}: linear, quadratic, and exponential (Figure \ref{fig:kernel}):
\begin{equation*}
    \kernel_\gamma (x) = \begin{cases}\frac{2}{\gamma^2}\left(\gamma-x\right) & \text{linear-type kernel}\\
      \frac{3}{2\gamma^3}\left(\gamma^2-x^2\right) & \text{quadratic-type kernel}\\
     \frac{e^{1/(x-\gamma)}}{\text{Ei}(-\gamma)+\gamma e^{-1/\gamma}} & \text{exponential-type kernel} ,
     \end{cases}
\end{equation*}
where $\text{Ei}$ is the exponential integral. For a given kernel type, we test several values of $\gamma$. For the PeMs dataset, we test $\gamma = \{0.004, 0.006, 0.1, 0.2, 0.3\}$. These values correspond to approximately $\{30, 45, 750, 1500, 2250\}$ meters, accounting for the scaling. For the drone-collected dataset, we test $\gamma = \{0.04, 0.06, 0.1, 0.2, 0.3\}$. These values correspond to approximately $\{32, 48, \\ 80, 160, 240\}$ meters. Before calibrating the model, we only keep data points that fall within $\left[0,1-\gamma\right]$ to ensure we only use existing data when estimating the convolution. For each combination of kernel type and $\gamma$, we perform the calibration process described in section \ref{sub:fd-cal}. We also present the two metrics, accuracy and coverage. The results of the  calibration on the corresponding datasets are shown in Table \ref{T:PeMS_result} and \ref{T:Drone_result}, respectively. 

\begin{table}[ht!]
\begin{center}
\caption{Free flow data from PeMS dataset. Highlighted rows are best-performing parameters.}
\label{T:PeMS_result}
\begin{tabular}{ |c|c|l|c|c|c| } 

\hline
Model & Kernel & $\gamma$ & $\kappa$ & \% Coverage ($\Sigma$) & \% Accuracy ($\varepsilon$)\\
\hline
\multirow{5}{*}{Nonlocal} & \multirow{5}{*}{exp} & \cellcolor{lightgray}0.004 (30 m)&\cellcolor{lightgray}0.6 &\cellcolor{lightgray}85 &\cellcolor{lightgray}81\\
& & 0.006 (45 m)& 0.4 & 80 & 79\\
& & 0.1 (750 m)& 0.4 & 79 & 78 \\
& & 0.2 (1500 m)& 0.3 & 79 & 78\\
& & 0.3 (2250 m)& 0.1 & 80 & 78 \\
\hline
\multirow{5}{*}{Nonlocal} & \multirow{5}{*}{linear} & 0.004 & 0.4 & 80 & 79\\
& & 0.006 & 0.4 & 80 & 79\\
& & 0.1 & 0.1 & 79 & 78 \\
& & 0.2 & 0.1 & 87 & 81\\
& & \cellcolor{lightgray}0.3 & \cellcolor{lightgray}0.0 & \cellcolor{lightgray}88 & \cellcolor{lightgray}81 \\
\hline
\multirow{5}{*}{Nonlocal} & \multirow{5}{*}{quadratic} & 0.004& 0.4 & 80 & 79\\
& &0.006 & 0.4 & 80 & 79\\
& & 0.1 & 0.0 & 79 & 79 \\
& & 0.2 & 0.0 & 87 & 81\\
& & \cellcolor{lightgray}0.3 & \cellcolor{lightgray}0.0 & \cellcolor{lightgray}90 & \cellcolor{lightgray}81 \\
\hline
\end{tabular}
\end{center}
\end{table}

\begin{table}[ht!]
\begin{center}
\caption{{Congestion data from drone dataset.  Highlighted rows are best-performing parameters.}}
\label{T:Drone_result}
\begin{tabular}{ |c|c|l|c|c|c| } 
\hline
Model & Kernel & $\gamma$ & $\kappa$ & \% Coverage ($\Sigma$) & \% Accuracy ($\varepsilon$)\\
\hline
\multirow{5}{*}{Nonlocal} & \multirow{5}{*}{exp} &\cellcolor{lightgray} 0.04 (32 m)
&\cellcolor{lightgray} 0.6 &\cellcolor{lightgray} 74 &\cellcolor{lightgray} 71\\
& &\cellcolor{lightgray} 0.06 (48 m)&\cellcolor{lightgray} 0.3 &\cellcolor{lightgray} 75 &\cellcolor{lightgray} 70 \\
& & 0.1 (80 m)& 0.2 & 73 & 68\\
& & 0.2 (160 m)& 0.3 & 69 & 39 \\
& & 0.3 (240 m)& 0.1 & 67 & 30 \\
\hline
\multirow{5}{*}{Nonlocal} & \multirow{5}{*}{linear} & \cellcolor{lightgray}0.04& \cellcolor{lightgray}0.0 & \cellcolor{lightgray}82 &\cellcolor{lightgray} 51\\
& &0.06& 0.5 & 67 & 38 \\
& &0.1& 0.9 & 55 & 24 \\
& &0.2& 0.1 & 35 & 14 \\
& &0.3& 1.0 & 32 & 8 \\
\hline
\multirow{5}{*}{Nonlocal} & \multirow{5}{*}{quadratic} & \cellcolor{lightgray} 0.04& \cellcolor{lightgray}0.0 & \cellcolor{lightgray}81 & \cellcolor{lightgray}48 \\
& &0.06& 0.5 & 63 & 34 \\
& &0.1& 0.7 & 54 & 22 \\
& &0.2& 0.3 & 28 & 9\\
& &0.3& 1.0 & 33 & 8 \\
\hline
\end{tabular}
\end{center}
\end{table}

\subsection{Results and interpretation}

Table \ref{T:PeMS_result} shows the result of calibrating the FD with the free flow data from PeMS. In this case, we consider the polynomial spline (see Section \ref{S:polynomial_aprroximation}) for the velocity profile. In the free flow scenario, under the exponential kernel, the model performs optimally for $\gamma = 0.004$ which corresponds to $30$ m of look-ahead information. An important observation is that for both linear and quadratic kernels in the free-flow scenario, the optimal calibration is associated with a physically unrealistic value of $\gamma = 0.3$ which corresponds to $2250$ m. Drivers are unlikely to be able to see long distances ahead. In these cases the calibrated $\kappa = 0$ refers to the nonlocal LWR model with no diffusion term, i.e. $\hat \rho = \rho$ in \eqref{E:fictitious_density}. Furthermore, more realistic values of $\gamma$ (i.e. 0.004 or 0.006) for the linear and quadratic kernels yield results that are not significantly worse than that for $\gamma = 0.3$. In a sense, while the data suggests multiple plausible fits, only a few are physically reasonable for interpreting traffic behavior. This means, in free flow scenario, the proposed model should be considered as best and most physically admissible model. 

Table \ref{T:Drone_result} shows the result of the accuracy and calibration in congestion scenario using the drone dataset. In this scenario, all kernels have the best fit when using $\gamma = 0.04$ (32 m). For the exponential kernel, $\gamma=0.06$ (48 m) performs similarly. For $\gamma = 0.04$, we find that $\{K_1=0.5, K_2=1.2, K_3=8.5\}$ and for $\gamma=0.06, \{K_1=0.9, K_2=1.0, K_3=8.6\}$. An important observation is that in the congestion scenario, in contrast to the free flow in Table \ref{T:PeMS_result}, larger $\gamma$ values, i.e. more look-ahead information, do not improve the accuracy of the model. This suggests that in congestion, drivers will consider immediate information much more strongly.

\begin{figure} [ht!]
    \centering
    \includegraphics[width=0.45\linewidth]{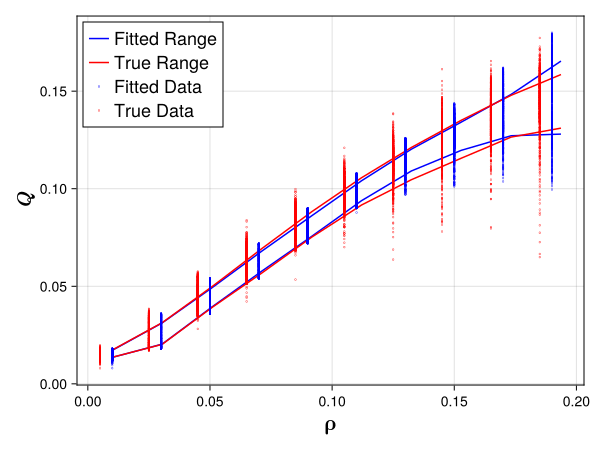}
    \includegraphics[width=0.45\linewidth]{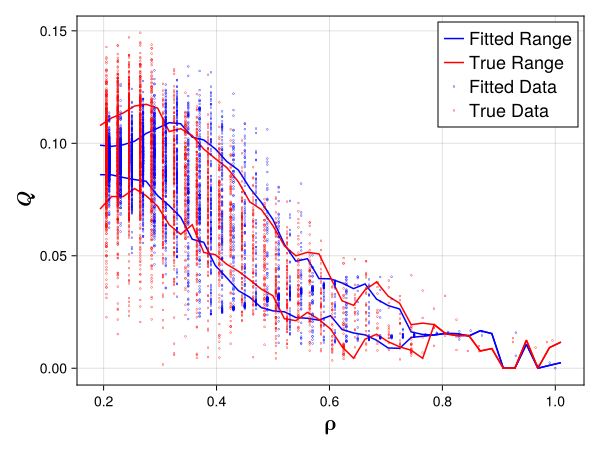}
    \caption{Calibrated fundamental diagram for (left) free flow condition using exponential kernel, $\gamma=0.004$ and $\kappa=0.6$: $\sum= 85$, $\epsilon=81$. (Right) Congested conditions using exponential kernel, $\gamma=0.04$, and $\kappa=0.6$: $\sum= 74$, $\epsilon=71$}
    \label{fig:calibrated_FD}
\end{figure}

\begin{figure}[ht!]
    \centering
    \includegraphics[width=0.45\linewidth]{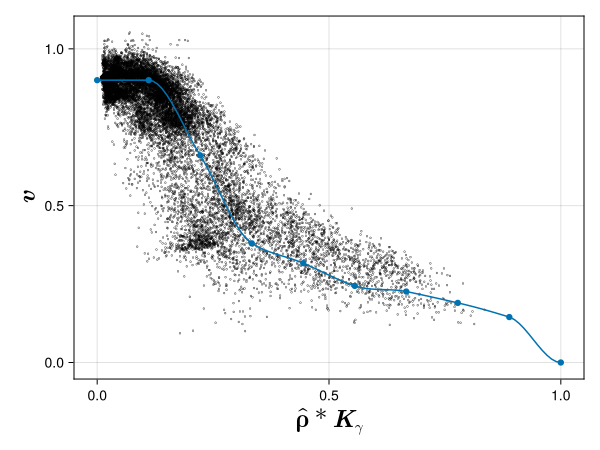}
    \includegraphics[width=0.45\linewidth]{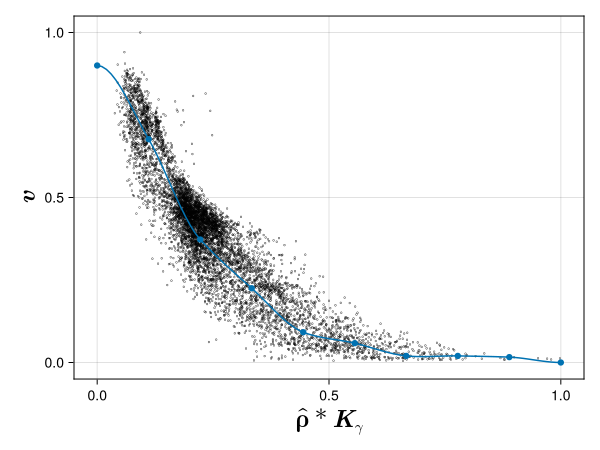}
    \caption{Calibrated polynomial velocity profile for (left) free flow condition using exponential kernel, $\gamma=0.004$, and $\kappa=0.6$. (Right) Congested conditions using exponential kernel, $\gamma=0.04$, and $\kappa=0.6$.}
    \label{fig:polynomial_velociy}
\end{figure}

Among both sets of results, the best parameter sets admit a non-zero value for the diffusive coefficient, $\kappa=0.6$. In other words, the proposed model provides more optimal fit than the previously studied nonlocal LWR (i.e. $\kappa=0.0$ is non-optimal), particularly for FD data. Additionally, the optimal physical distances, kernels, and coefficient $\kappa$ are similar for both sets of data, exponential kernel with $\gamma=30$ m and $\kappa=0.6$. These results suggest that drivers consider the same look-ahead information regardless of the traffic flow regime, and, potentially, independent of road characteristics. Figure \ref{fig:calibrated_FD} shows the optimal result for the free flow (left) and the congestion (right) data. Figure \ref{fig:polynomial_velociy} shows the calibrated velocity functions for each dataset.

\begin{figure}[ht!]
    \centering
    \includegraphics[width=0.45\linewidth]{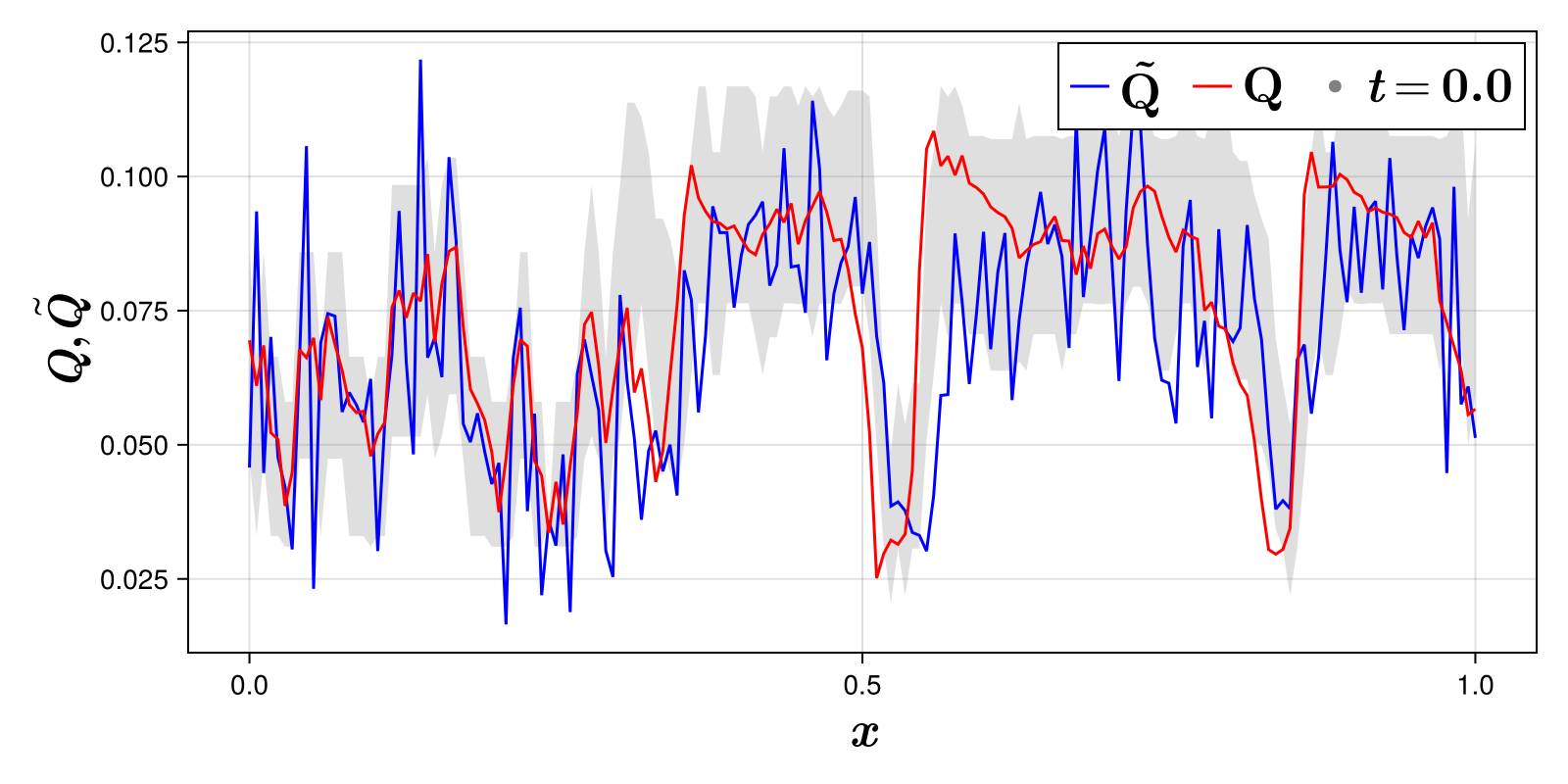}
    \includegraphics[width=0.45\linewidth]{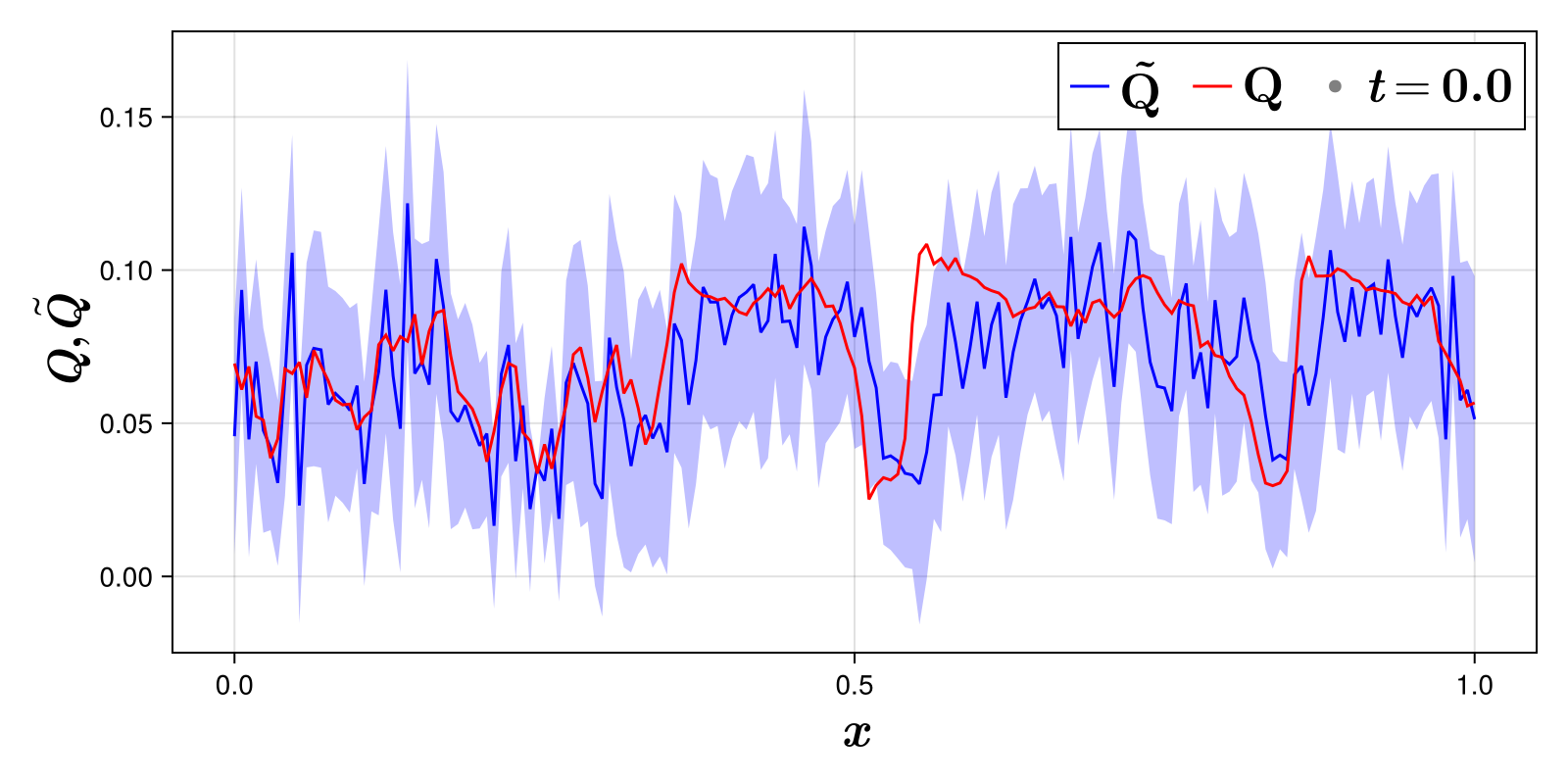}
    \includegraphics[width=0.45\linewidth]{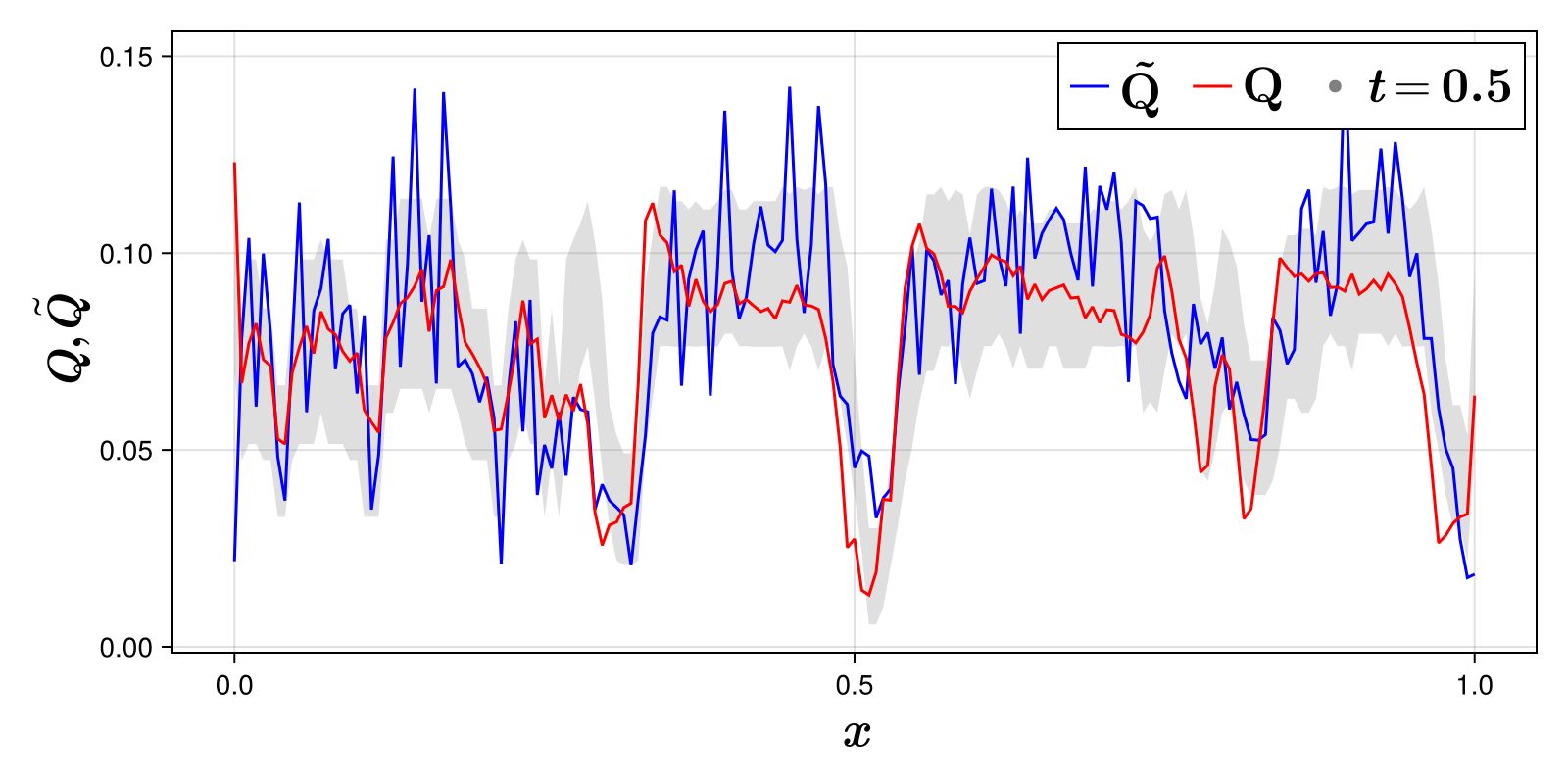}
    \includegraphics[width=0.45\linewidth]{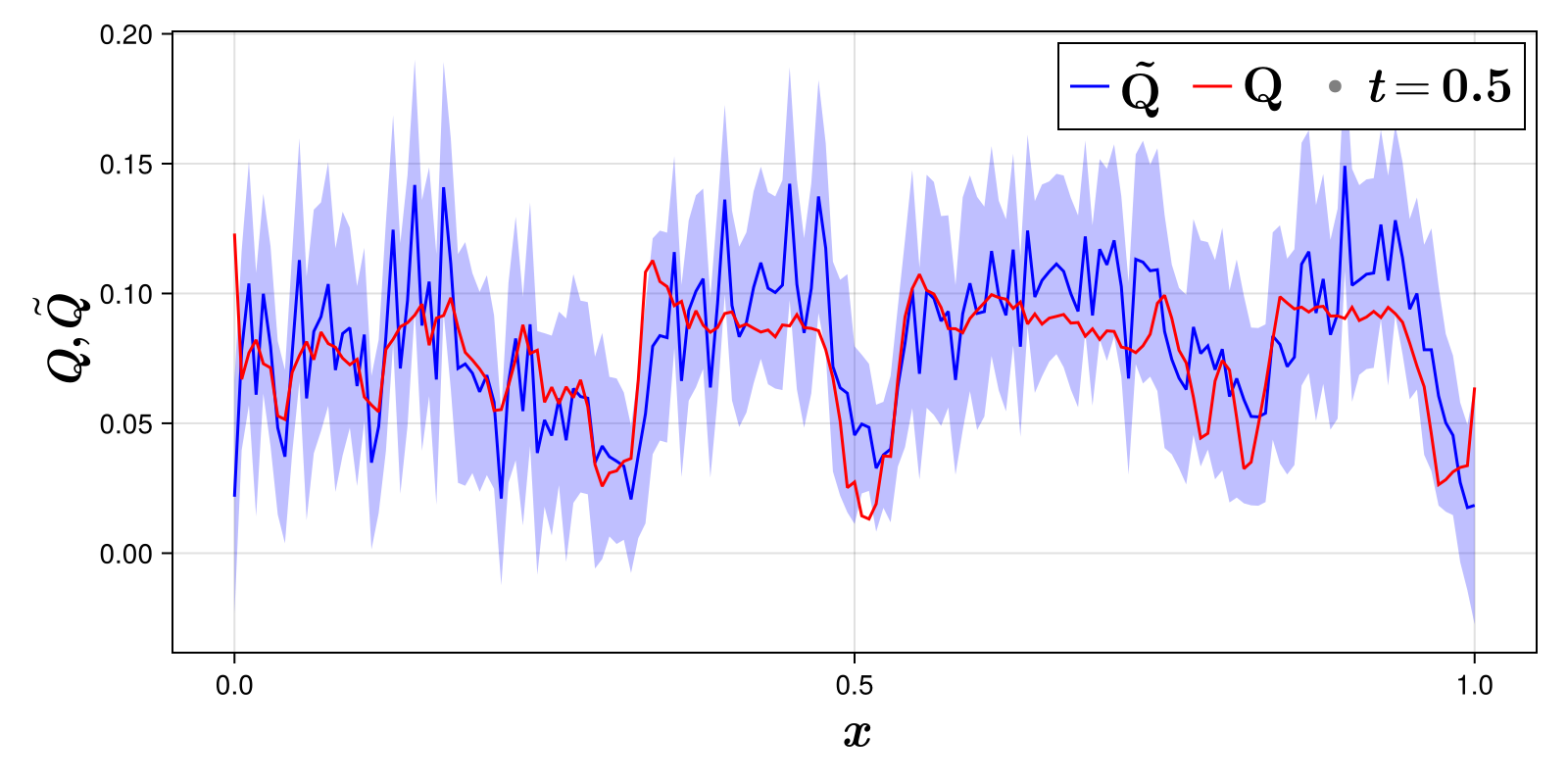}

    \includegraphics[width=0.45\linewidth]{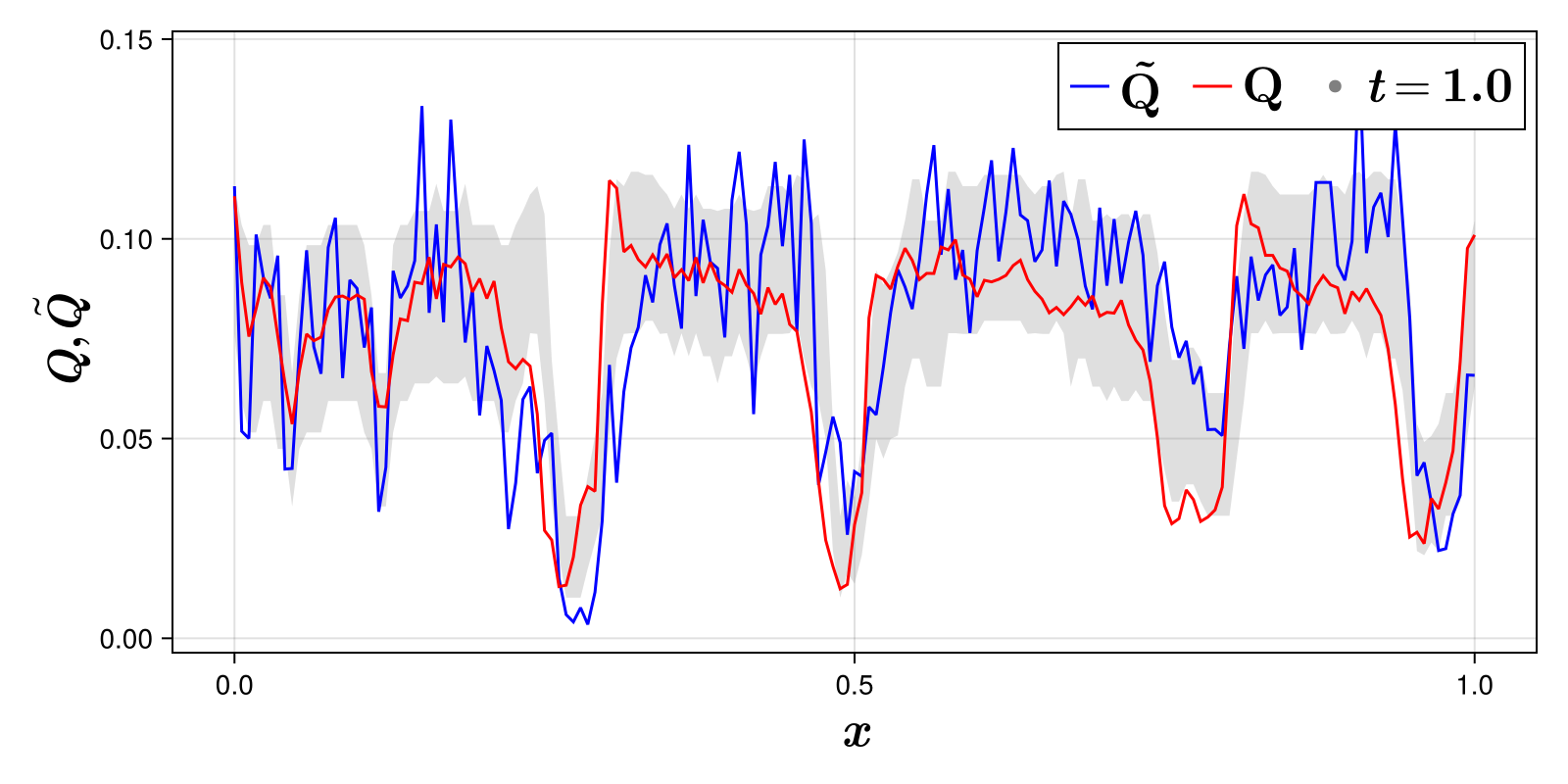}
    \includegraphics[width=0.45\linewidth]{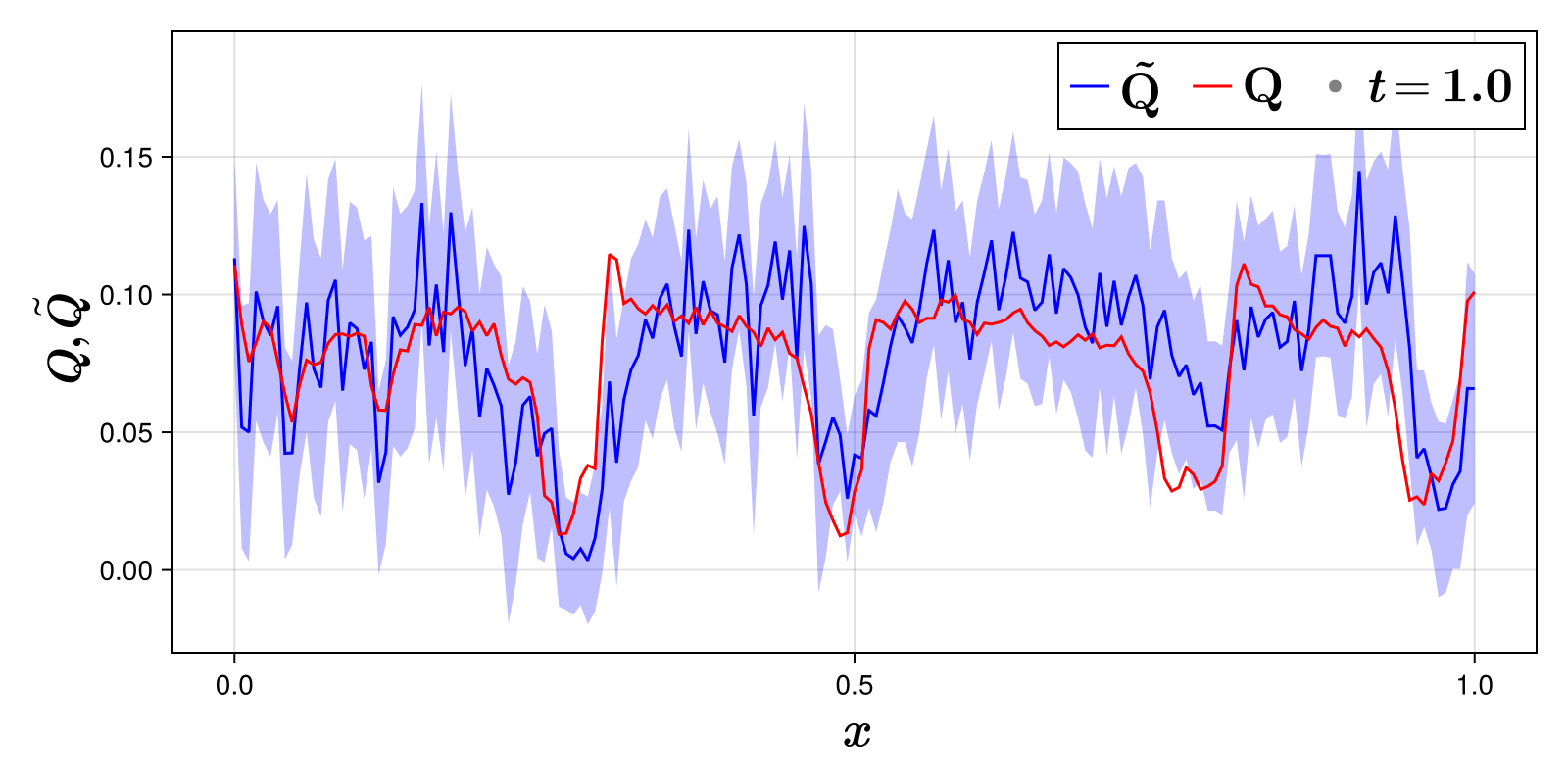}

    \includegraphics[width=0.45\linewidth]{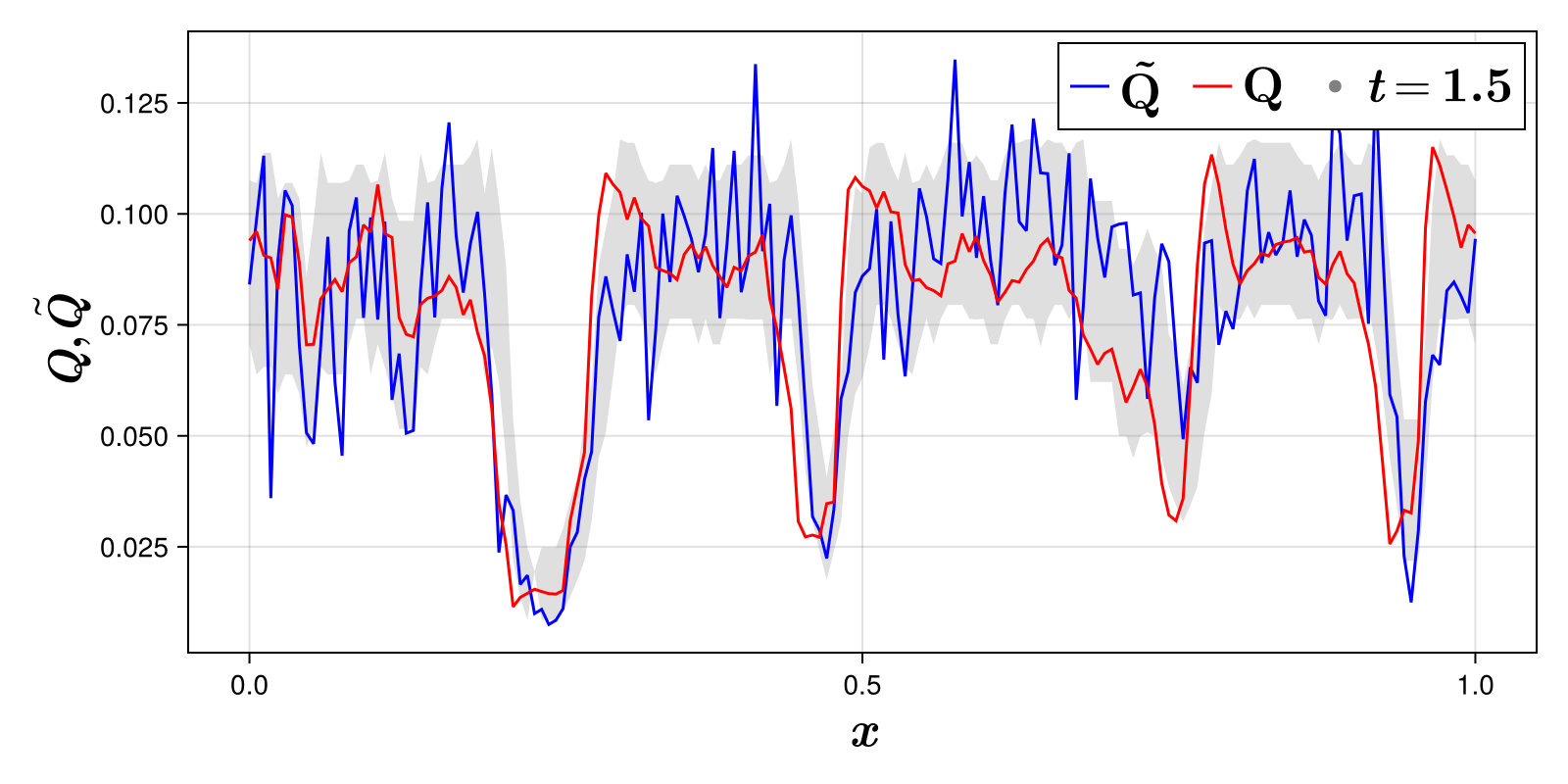}
    \includegraphics[width=0.45\linewidth]{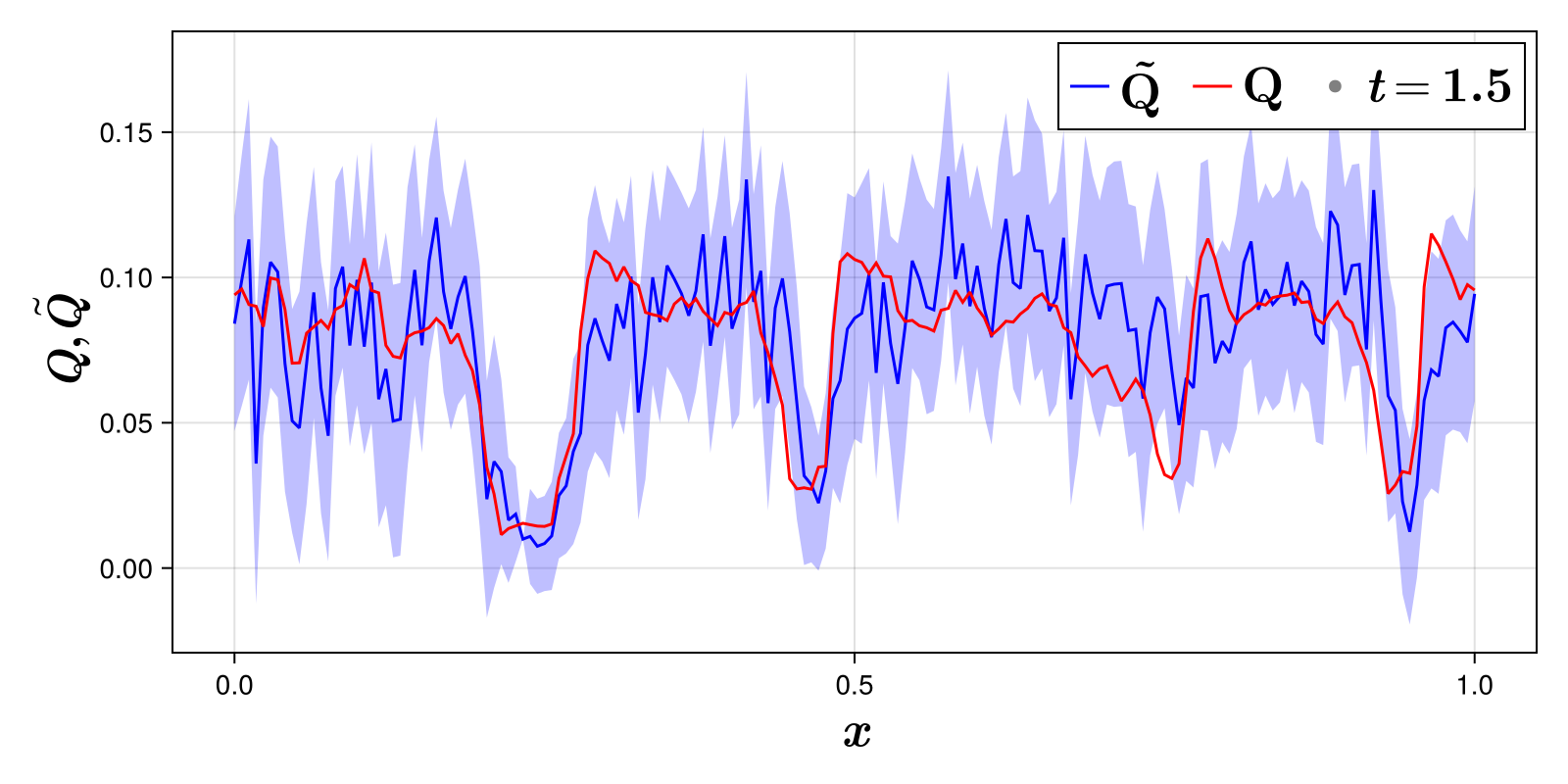}

    \includegraphics[width=0.45\linewidth]{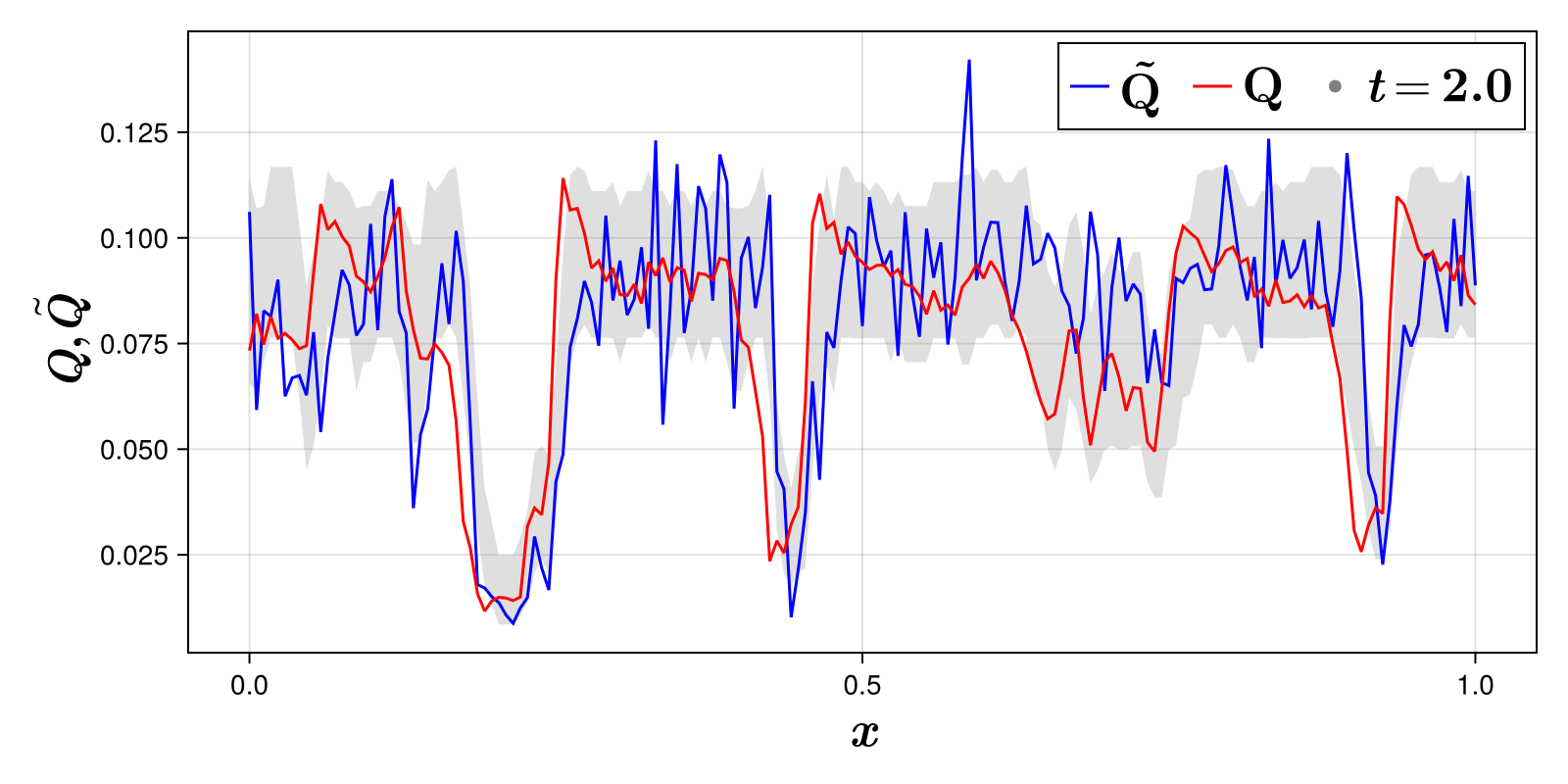}
    \includegraphics[width=0.45\linewidth]{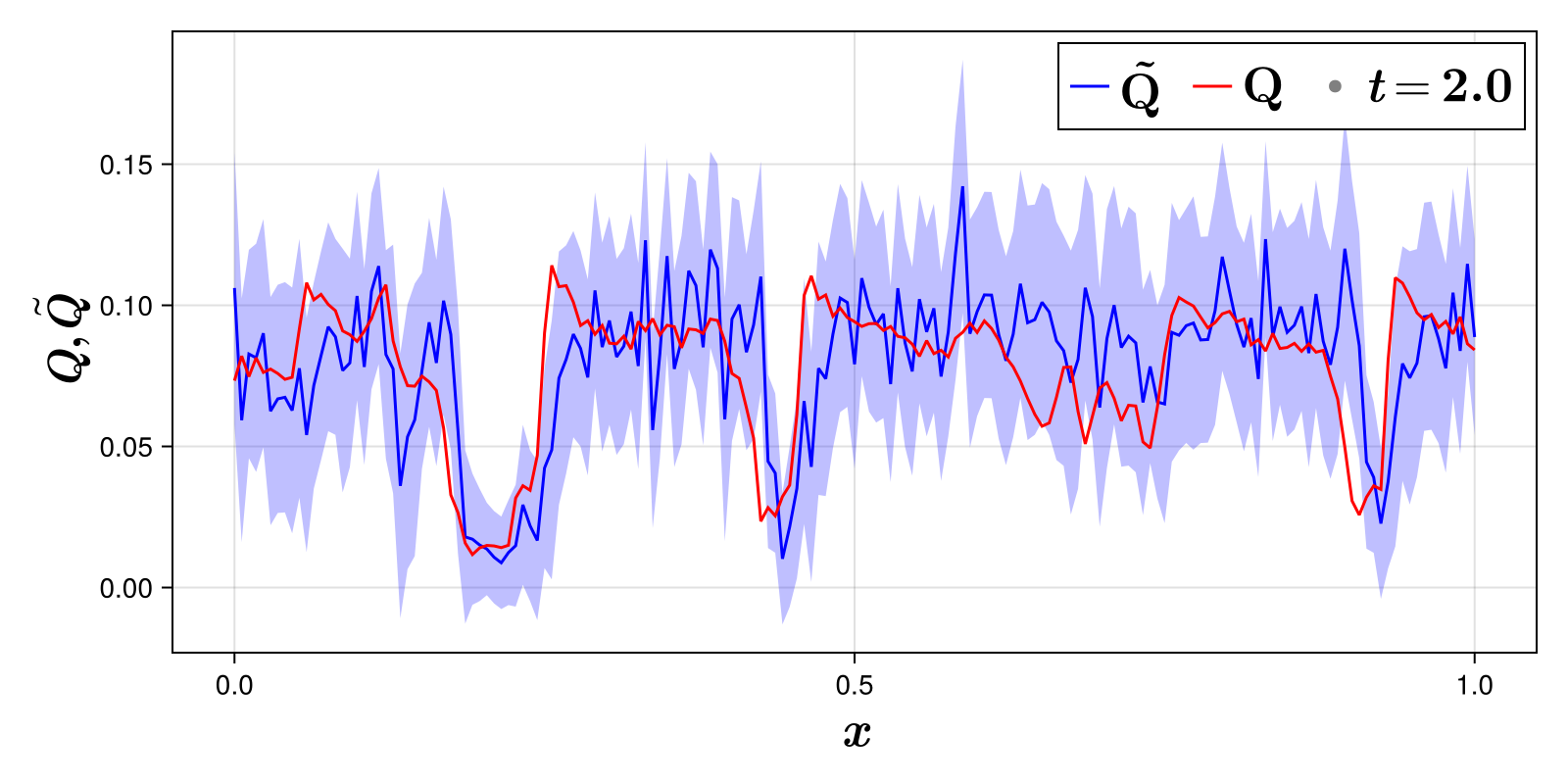}
    
    \caption{Comparison between actual and predicted flow at different times. (Left) gray bands represent $\mu_{0,m} \pm S_{0,m}$, (right) blue bands represent $\tilde{Q}\pm 2S_{0,m}$}
    \label{fig:comparison_Q}
\end{figure}

We also compare the predicted flow profiles, $Q= f(\tilde{\rho}(t,x);\mathbf{P}^*)$, with the measured flow profiles $\tilde Q(t,x)$ for arbitrary times (\ref{fig:comparison_Q}). It should be noted that while our calibration method does not guarantee that the error between the ground-truth flow measurement and estimated flow for each measurement will be close (unlike a least-squares based method), Figure \ref{fig:comparison_Q} shows that the calibrated flows still closely reproduce them. Along with the flows, we plot intervals showing the prediction range within 1, 1.5, and 2 standard deviations. The subfigures on the left in Figure \ref{fig:comparison_Q} show the intervals defined by $\mu_{0,m} \pm S_{0,m}$, where $m$ identifies the range that contains $\tilde \rho (t,x)$. These intervals show the range of expected flows given by the empirical data. The subfigures on the right in Figure \ref{fig:comparison_Q} show the intervals defined by $\tilde Q(t,x) \pm 2S_{0,m}$, with $m$ defined similarly, show an estimated error margin. At a glance, the predicted flows typically fit the shape measured flows and fall within the plotted bands. The largest discrepancies are often when the flow transitions from low to high within a short distance. These samples correspond to when the density transitions from super-critical to near-critical (i.e. when the derivative of density is a high-magnitude negative value). Table \ref{comparison_q} enumerates the percentage of predicted points that fall outside each type of band, for several time samples. We can see that the predicted flow lies within a bound of just 1 standard deviation of average measured flow over 80\% of the time, up to over 98\% for 2 standard deviations. This suggests that with this calibration method, the predicted flows indeed will closely correspond to expected flows, which is desirable for a macroscopic model. When compared to the bands which represent error margins, the predicted flows fall within 2 standard deviations up to about 90\% of the time. These results suggest that even for point-wise prediction (i.e. in a non-macroscopic sense), this calibration method can produce accurate results. 

\begin{table}[ht!]
    \centering
    \caption{Accuracy of Predicted Flow}
    \label{comparison_q}
    \begin{tabular}{|c|c|c|c|c|c|c|}
        \hline
        \textbf{time} & \multicolumn{6}{c|}{\textbf{\% samples outside the band}} \\
        \hline
        $[hr]$ & $\mu_{0,m} \pm S_{0,m}$ & $\mu_{0,m} \pm 1.5S_{0,m}$ & $\mu_{0,m} \pm 2S_{0,m}$ & $\tilde{Q} \pm S_{0,m}$ & $\tilde{Q} \pm 1.5S_{0,m}$ & $\tilde{Q} \pm 2S_{0,m}$ \\
        \hline
        0.0 & 16.18 & 4.05 & 1.73 & 34.1 & 16.18 & 7.51 \\
        0.5 & 17.34 & 8.09 & 2.31 & 35.84 & 18.5 & 7.51 \\
        1.0 & 18.5 & 8.09 & 2.89 & 42.2 & 25.43 & 11.56 \\
        1.5 & 19.65 & 10.98 & 4.05 & 41.62 & 20.23 & 10.4 \\
        2.0 & 17.92 & 8.67 & 2.89 & 37.57 & 19.08 & 10.4 \\
        \hline
    \end{tabular}
\end{table}

\section{Results for solution-based calibration}

\begin{figure}[ht!]
    \centering
    \includegraphics[width=0.3\linewidth]{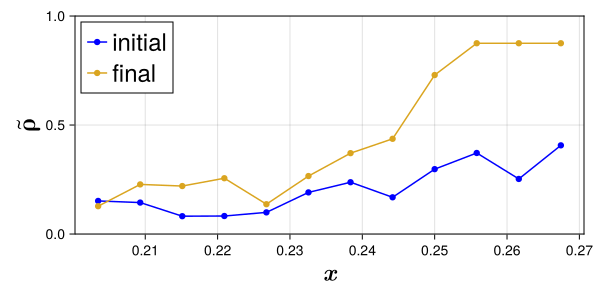}
    \includegraphics[width=0.3\linewidth]{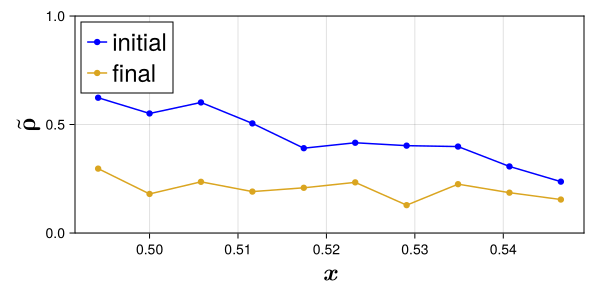}
    \includegraphics[width=0.3\linewidth]{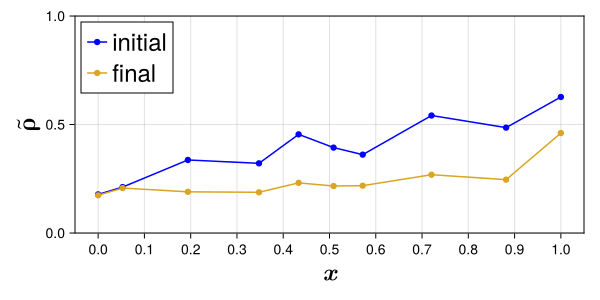}
    \caption{Measured density profiles chosen for the performance calibration experiments. (Left) Congestion from drone data, (middle) dissipation from drone data, (right) dissipation from PeMS.}
    \label{fig: congestion_dissipation}
\end{figure}

In this section, we describe the calibration of the model with respect to the accuracy of the PDE solution to real world traffic situations. We assess the model's performance in a congestion and dissipation from the drone-based dataset. We also use a dissipation scenario from the PeMS dataset as an additional experiment. These scenarios were chosen to test the model's performance in super-critical conditions (a known weakness of LWR).  
For the congestion scenario, we selected the section from 0.2 to 0.27 in the Figure \ref{fig:x_t_plot_heatmap}. Additionally, we defined a time frame of 15 minutes, approximately from 4:55 PM to 5:10 PM, to analyze the model's performance during peak congestion.
For the dissipation scenario, the entire section of the PeMS data was utilized. The initial time was set at 8:30 AM, and the final time was 10:00 AM. This time frame captures the gradual reduction in traffic density, allowing us to analyze the model’s ability to represent the dissipation phase. 
For the dissipation scenario, another section was selected from the drone dataset. The location ranges from 0.49 to 0.55, as shown in Figure \ref{fig:x_t_plot_heatmap}. The selected time frame is 1 hour, spanning from 4:55 PM to 5:55 PM. The chosen section is shown in Figure \ref{fig: congestion_dissipation}.

For each scenario, we calibrate the parameters of the model using the method described in Section \ref{sec:solution-based}. We also perform the calibration on LWR (i.e. $\gamma=\kappa=0$) and the nonlinear diffusively-corrected LWR \cite{campos2021saturated, colombo2012class} for comparison. The nonlinear diffusively-corrected LWR (presented with the alias, ``Phi," in following figures and tables) is defined as follows:
\begin{equation}\label{eq:nonlinear}
     \partial_t \rho + \partial_x [\rho U(\rho)-\kappa D(\rho)\Psi(\partial_x \rho)] = 0.
\end{equation}
In Table \ref{T:performance}, we present our optimal results and parameters for the model using Newell's velocity function. For each model, we also present the mean squared residual (MSR):
\begin{equation*}
    MSR = \frac{1}{N\cdot T}\sum_i^T\sum_j^N (\tilde\rho(i\Delta t,j\Delta x)-\rho(i\Delta t,j\Delta x))^2
\end{equation*}

\if{
\begin{figure}[ht!]
    \centering
    \includegraphics[width=0.31\linewidth]{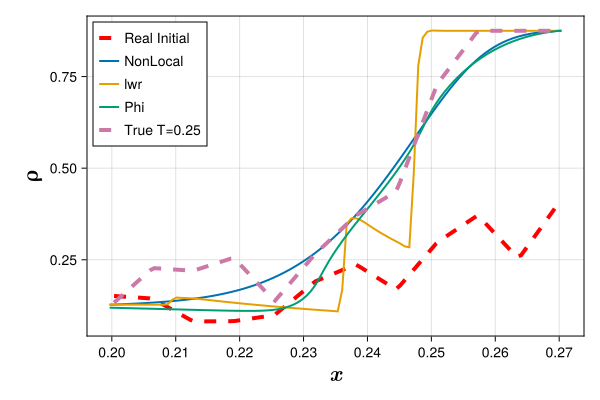}
    \includegraphics[width=0.31\linewidth]{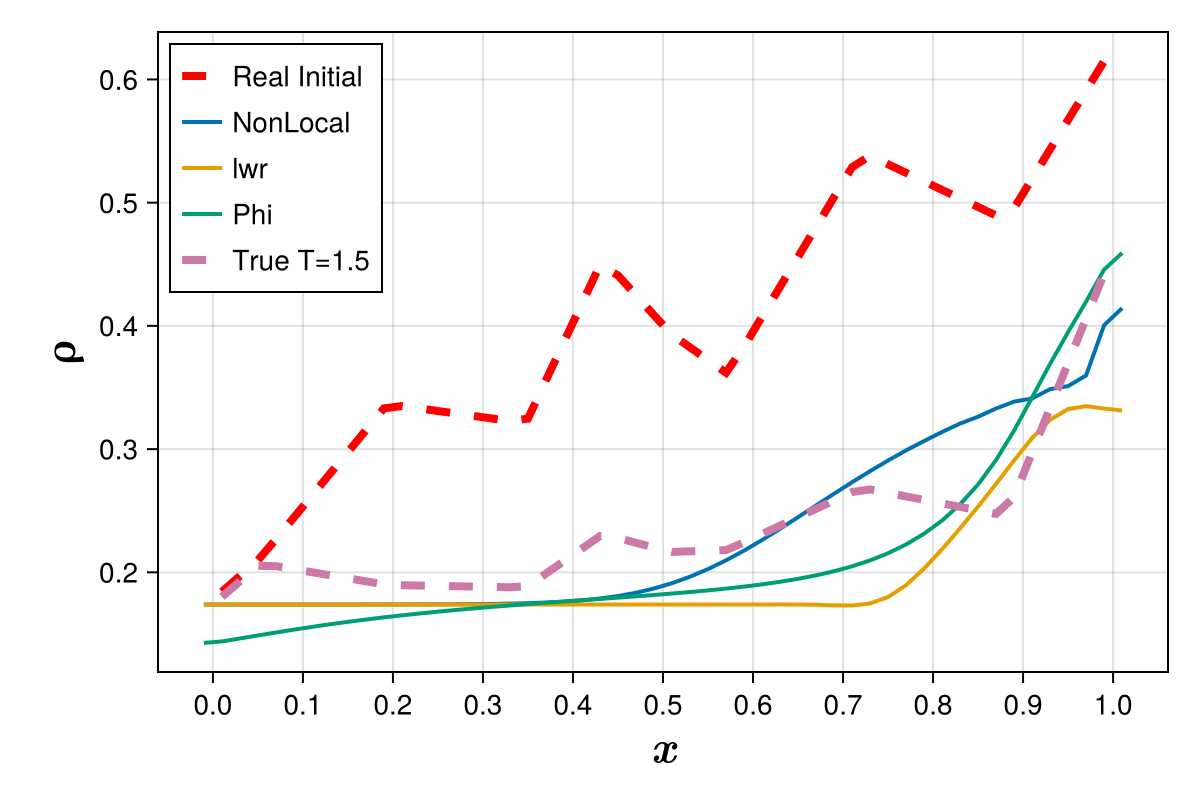}
    \includegraphics[width=0.31\linewidth]{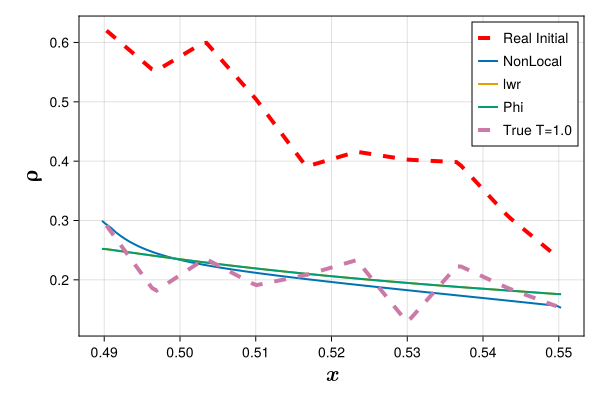}
    \caption{Performance using different models [Left-Congestion(drone), Middle-Dissipation(PeMS), Right-Dissipation(Drone)]}
    \label{fig: performance}
\end{figure}}
\fi

\begin{figure}[ht!]
    \centering
    \includegraphics[width=0.32\linewidth]{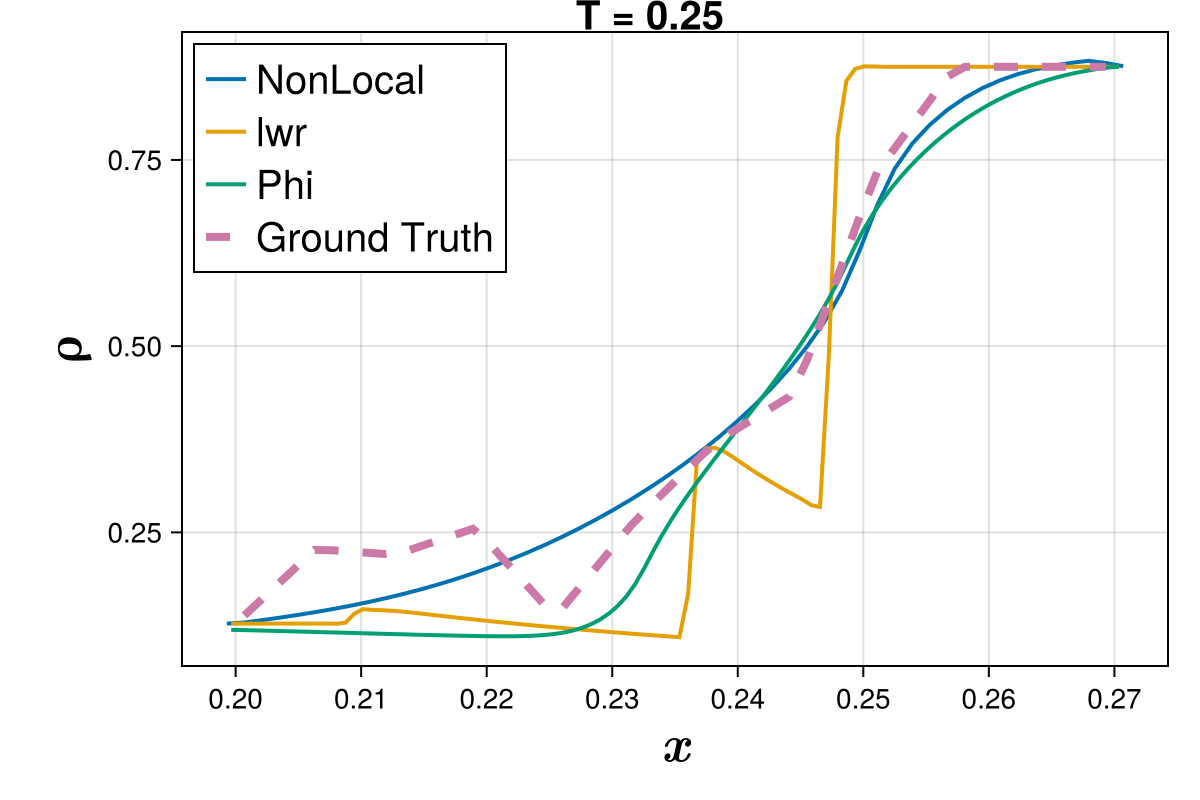}
    \includegraphics[width=0.32\linewidth]{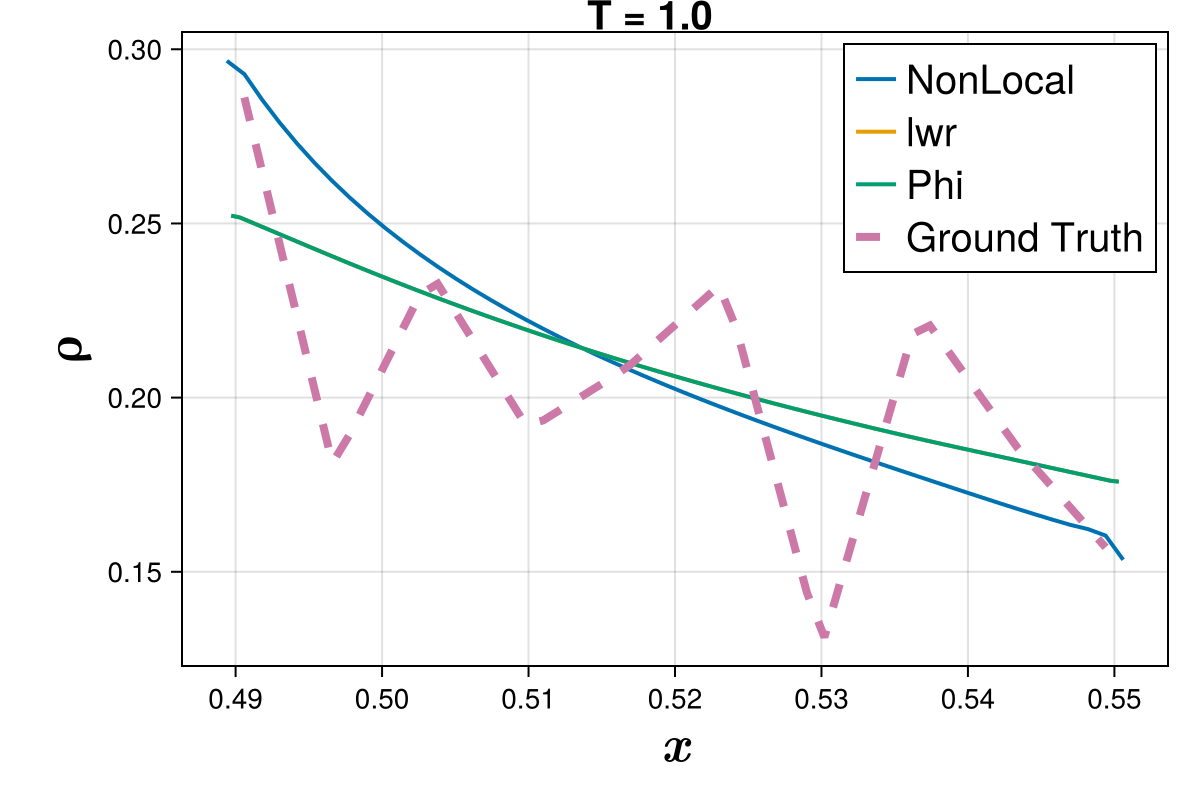}
    \includegraphics[width=0.32\linewidth]{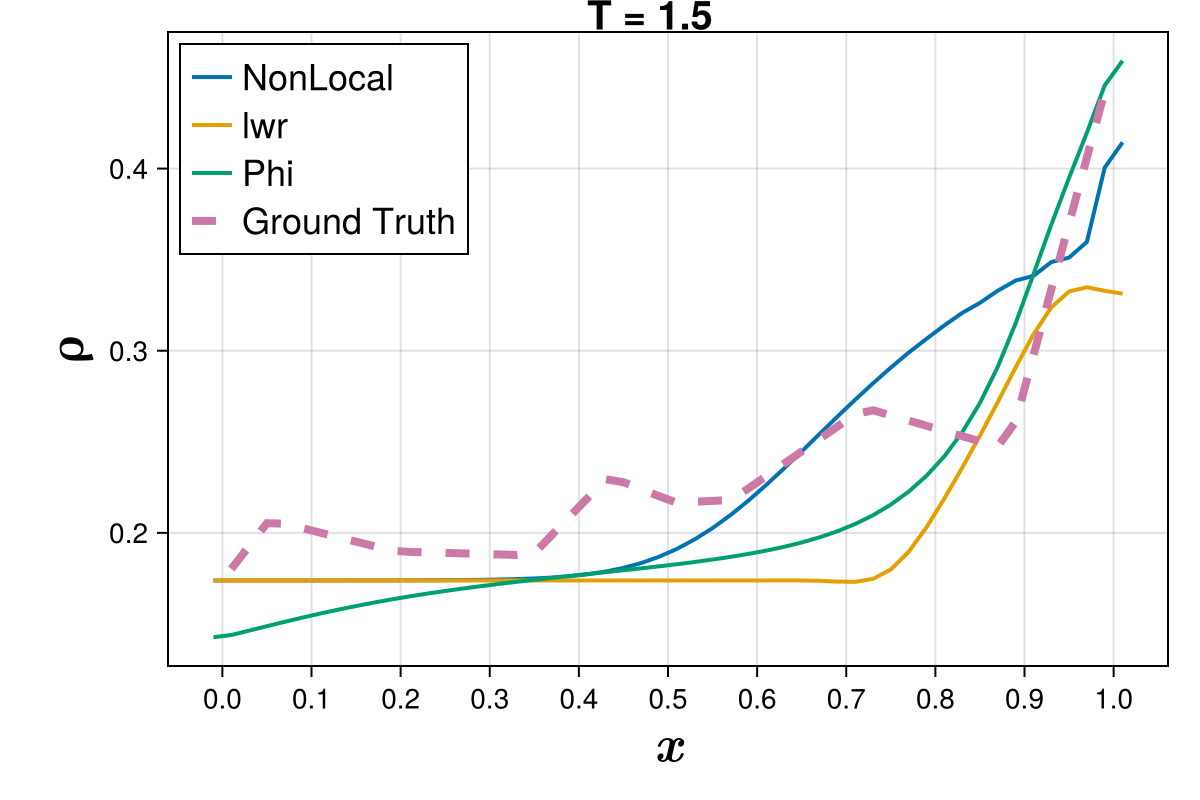}
    
    \caption{Optimal performance (final time comparison) of the proposed nonlocal model, LWR, and the nonlinear diffusively-corrected LWR for (left) congestion, (middle) dissipation-drone, (right) dissipation-drone experiments.}
    \label{fig:gamma_04}
\end{figure}

\begin{table}[ht!]
    \centering
    \caption{Performance for Different Models }
    \label{T:performance}
    \begin{tabular}{|c|c|c|c|c|c|c|c|c|}
    \hline
    Dataset & Scenario & Kernel & Model & $\gamma$ &$\kappa$ & MSR & $v_{max}$ & $c$ \\
    \hline
    \multirow{4}{*}{Drone} &\multirow{4}{*}{Congestion} &  \multirow{4}{*}{exp}& Nonlocal & 0.003 & $0.3$ & $0.002009$ & $1.8$ & $0.1$ \\ 
    & & & & 0.01 & $0.1$ & $0.002794$ & $1.3$ & $0.1$\\
     & & & LWR & & $0.0$ & $0.009996$ & $1.4$& $0.1$\\
     & & & Phi & & $0.8$ & $0.005074$ & $1.0$ &$0.1$ \\
     \hline
     \multirow{3}{*}{Drone}& \multirow{3}{*}{Dissipation} & \multirow{3}{*}{exp} & Nonlocal& 0.003 & $0.0$& $0.000479$ & $1.8$ & $0.2$ \\
     & & & & 0.01 & $0.0$& $0.000725$ & $1.6$ & $0.1$ \\
     & & & LWR & & $0.0$ & $0.000638$ & $1.5$ & $0.1$\\
     & & & Phi & & $0.0$ & $0.000638$ & $1.5$ & $0.1$\\
     \hline
     \multirow{3}{*}{PeMS}& \multirow{3}{*}{Dissipation} & \multirow{3}{*}{exp} &
     Nonlocal & $0.004$ & $1.0$& $0.00124$ & $1.3$ & $1.0$ \\
     & & & & $0.2$ & $1.0$& $0.00115$ & $1.1$ & $0.9$ \\
     & & & LWR & & $0.0$ & $0.00224$ & $0.5$ & $1.9$\\
     & & & Phi & & $0.1$ & $0.00145$ & $0.7$ & $0.9$\\
     \hline
    \end{tabular}
\end{table}

\subsection{Results and interpretation}
\begin{figure}[ht!]
    \centering
    \includegraphics[width=0.32\linewidth]{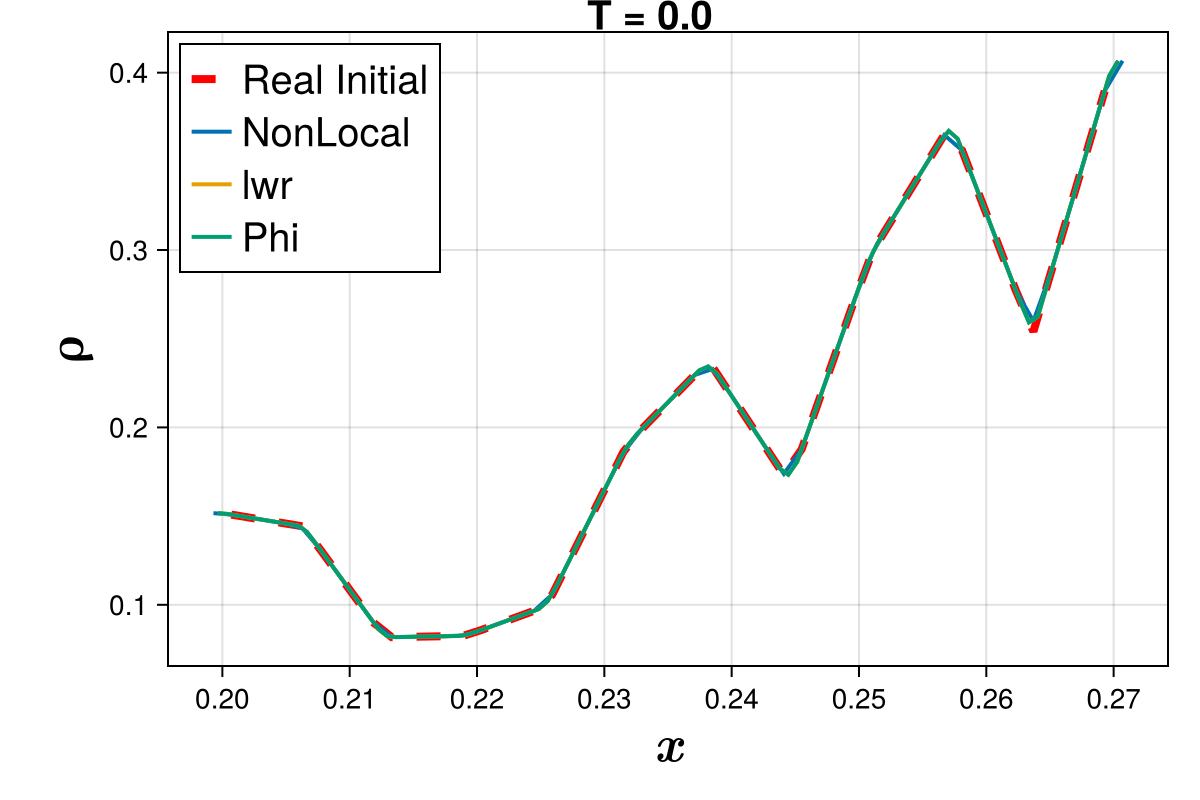}
    \includegraphics[width=0.32\linewidth]{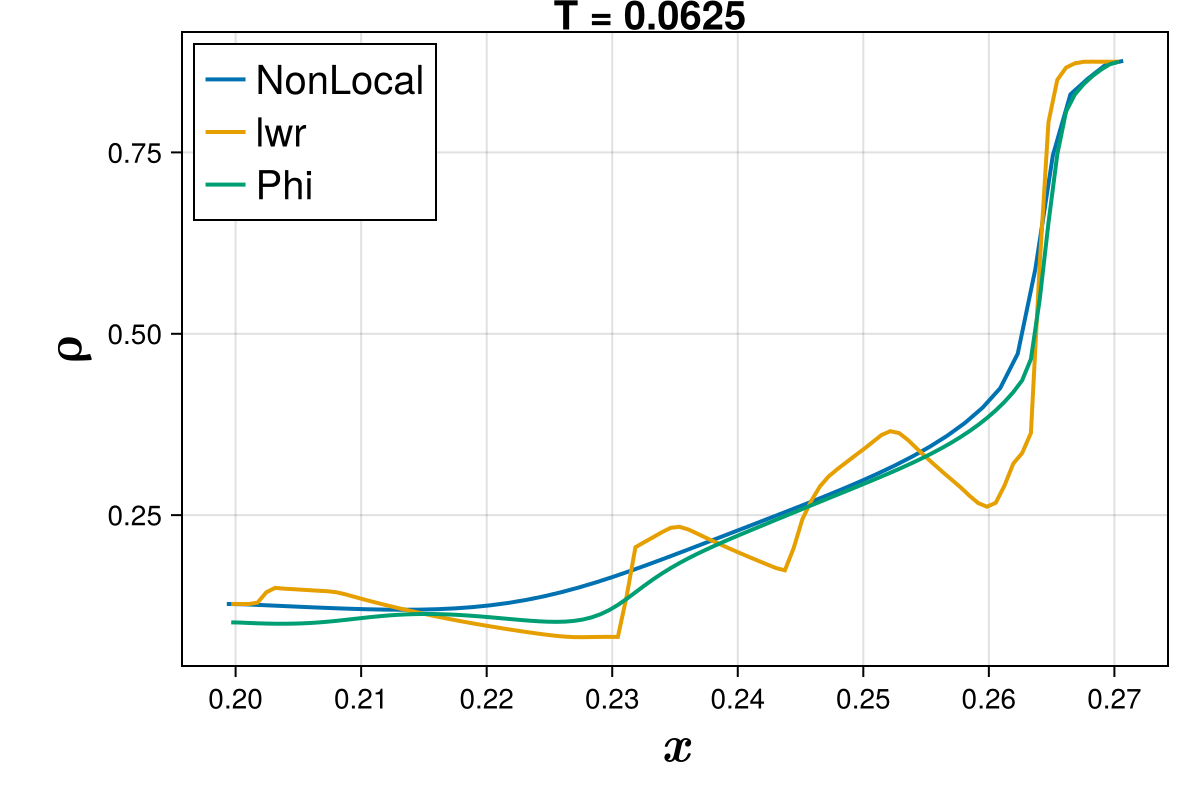}
    \includegraphics[width=0.32\linewidth]{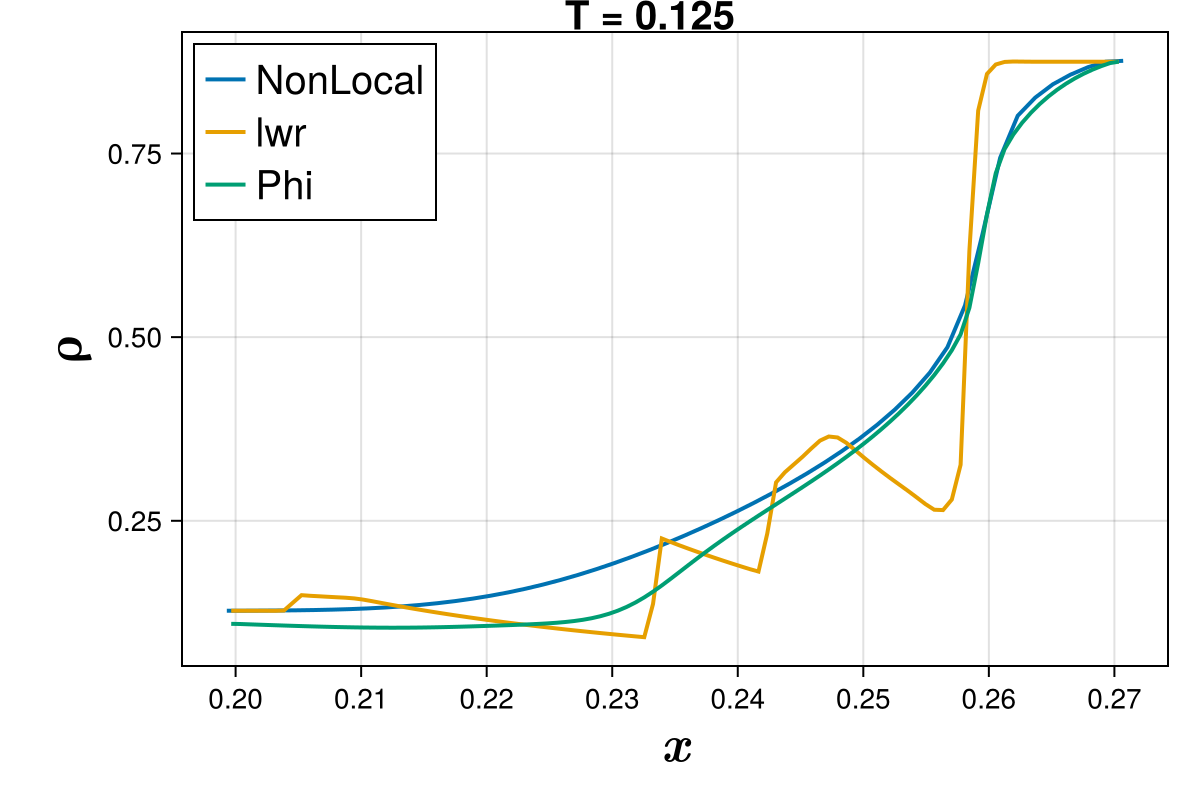}
    \includegraphics[width=0.32\linewidth]{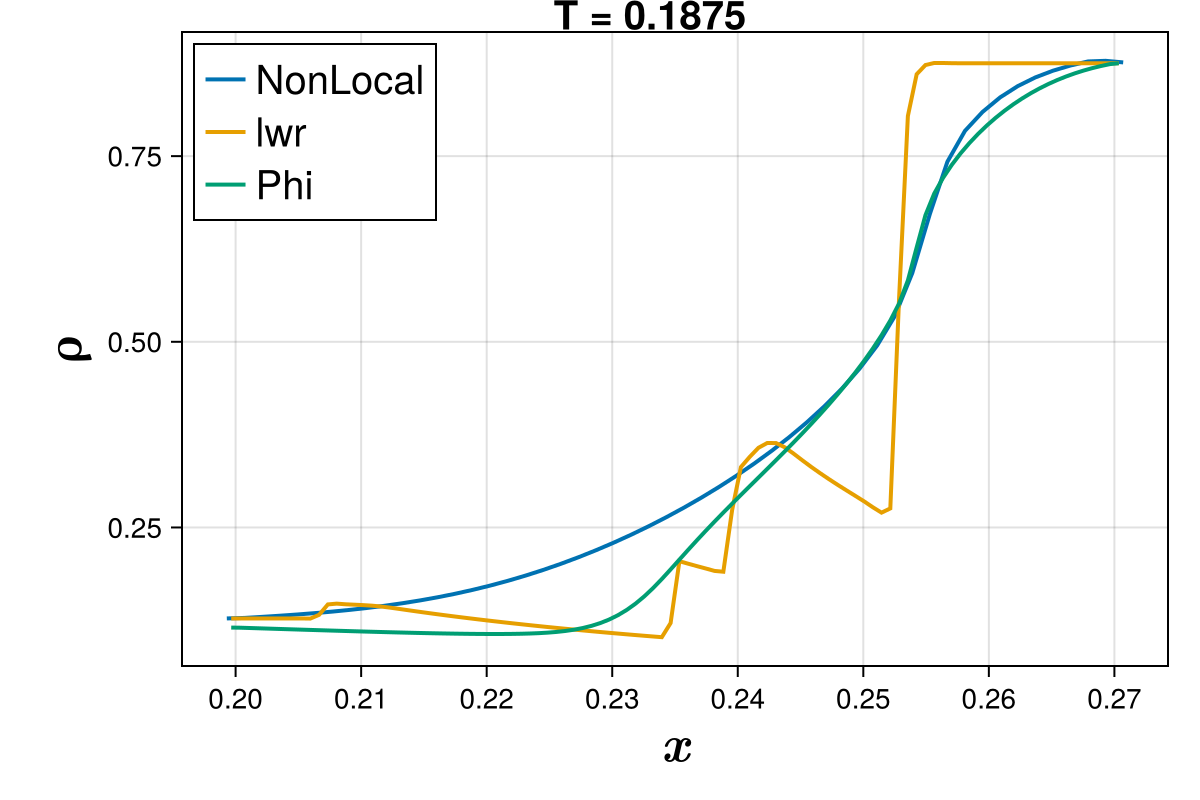}
    \includegraphics[width=0.32\linewidth]{Figures/Performance_final/congestion_gamma_04_q4.png}
    \caption{\small{Step-by-step progression of performance in congestion using the drone data. The frames correspond to simulation times of \( T = 0\), \(0.0625 \), \( 0.125 \), \( 0.1875 \), and \( 0.25 \) hr, where \( T \) represents different time steps. The ground truth is the final density value from the data at \( T_G = 0.25 \) hr.}}
\label{fig:congestion_evulotion}
\end{figure}

Figure \ref{fig:gamma_04} shows the optimal solution of the three models at the final time steps of their respective scenarios, along with the initial and final measured density profiles. Table \ref{T:performance} shows that the mean squared residual of the proposed model at $\gamma=0.003, 0.01$ ($2, 11$ m) is about 25\% of the error of LWR and half the error of the nonlinear diffusive model, \eqref{eq:nonlinear}, in congestion. Figure \ref{fig:congestion_evulotion} shows that the evolution of congestion from the nonlocal model can capture the same density profile of the measured data and matches the shape of the shockwave better than LWR. Compared to the nonlinear diffusive model \eqref{eq:nonlinear}, the nonlocal model also adheres closely in sub-critical flow. 

In the dissipation scenario using drone data, the error of the nonlocal model is 32\% lower than that of LWR and \eqref{eq:nonlinear} when $\gamma=0.003$, and similar when $\gamma=0.01$. In the dissipation scenario using PeMS data, the nonlocal model performs similarly to \eqref{eq:nonlinear} and LWR. In particular, we test the value for $\gamma$ found from the FD calibration. From both experiments, the three models are able to capture flow that transitions from supercritical to subcritical. Overall, these results show that the proposed model has a significantly greater ability to capture congested traffic phenomena. At worst, it performs similarly to previous models depending on the value chosen for $\gamma$. We note that in dissipation, traffic flow transitions into free flow, where the relationship between density and flow is significantly more linear (i.e. the left-hand-side of the FD has low variance), which can explain why the nonlocal and \eqref{eq:nonlinear} models degenerate to LWR in performance.
\begin{figure}[ht!]
    \centering
    \includegraphics[width=0.32\linewidth]{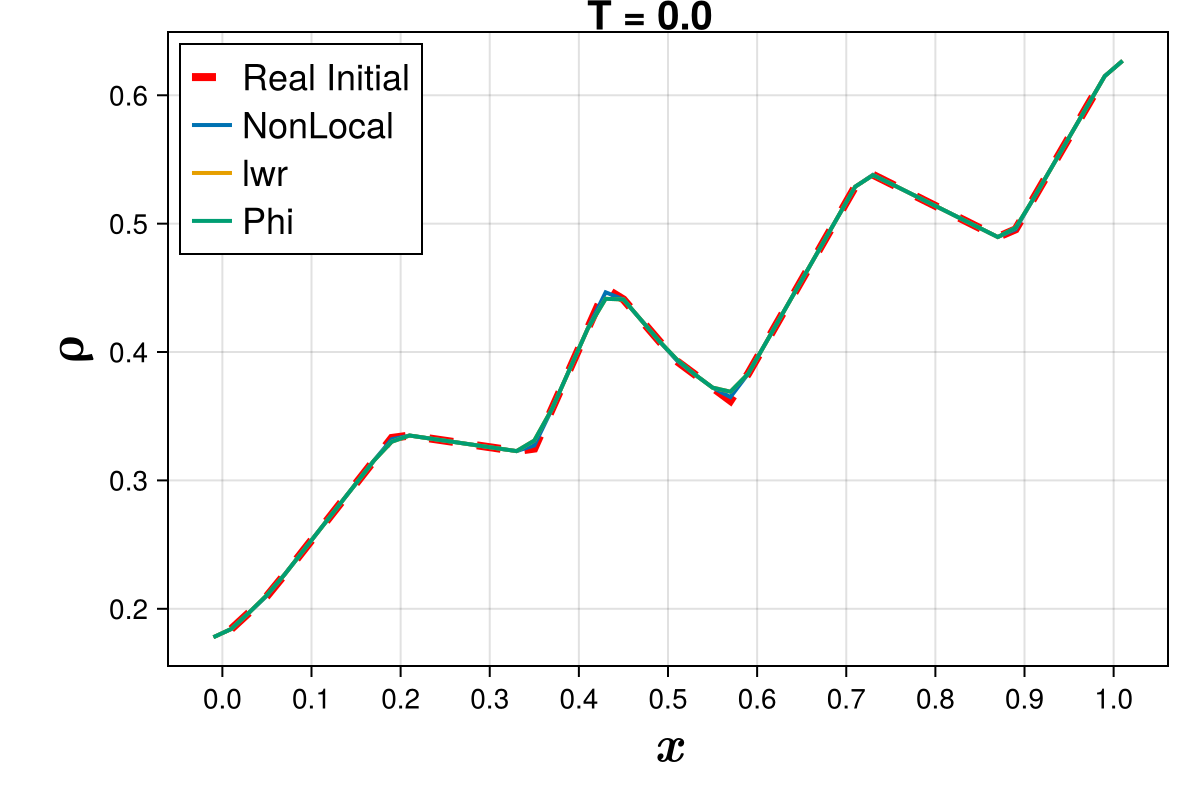}
    \includegraphics[width=0.32\linewidth]{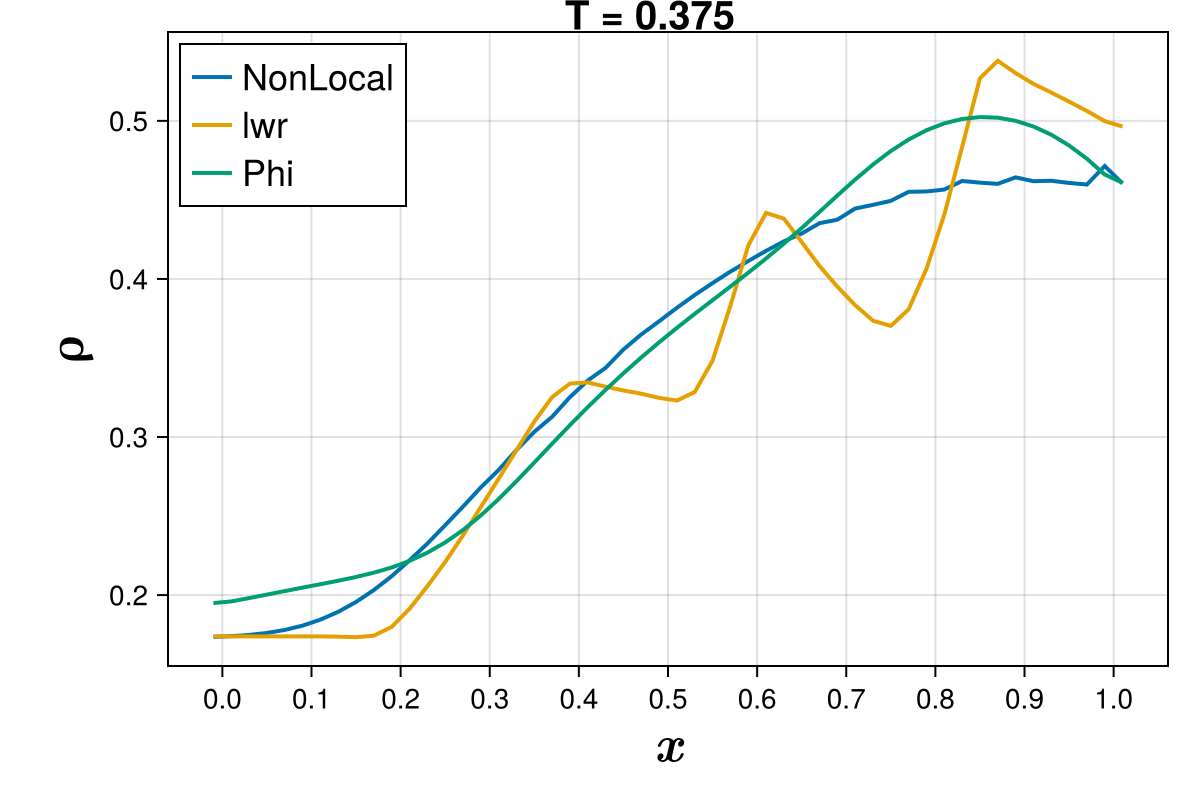}
    \includegraphics[width=0.32\linewidth]{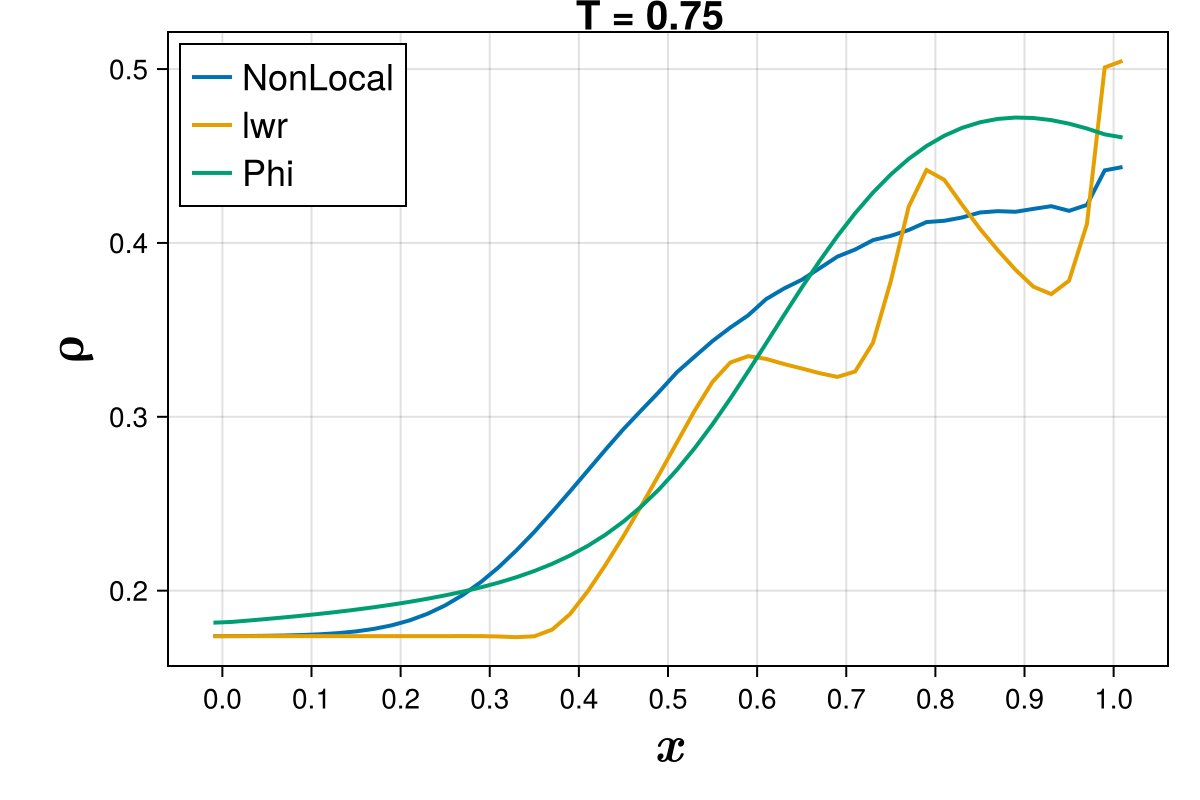}
    \includegraphics[width=0.32\linewidth]{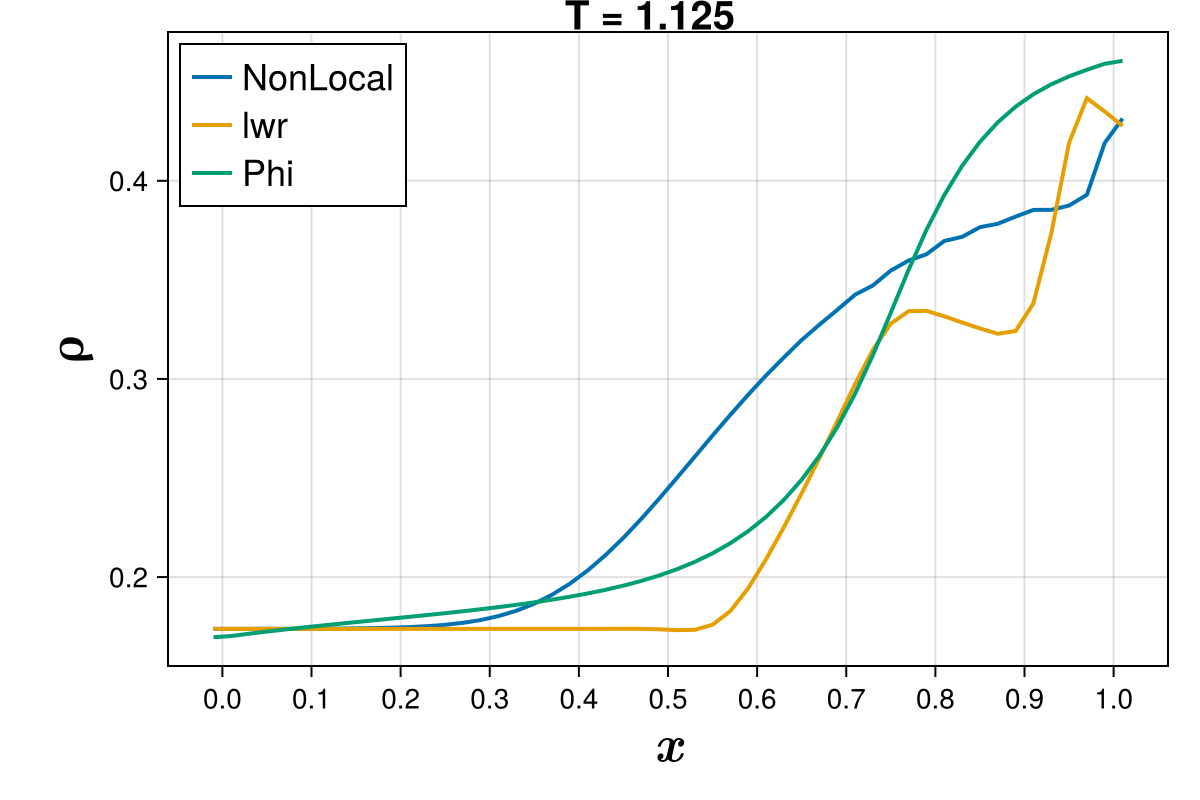}
    \includegraphics[width=0.32\linewidth]{Figures/Performance_final/dissipation_PeMS_gamma02_q4.png}
    \caption{\small{
    Step-by-step progression of performance in dissipation using the PeMS data. The frames correspond to simulation times of \( T =0\), \(0.375 \), \( 0.75 \), \( 1.125 \), and \( 1.5 \) hr, where \( T \) represents different time steps. The ground truth is the final density value from the data at \( T_G = 1.5 \) hr.}}

    \label{fig:dissipation_pems_evulotion}
\end{figure}
\begin{figure}[ht!]
    \centering
    \includegraphics[width=0.32\linewidth]{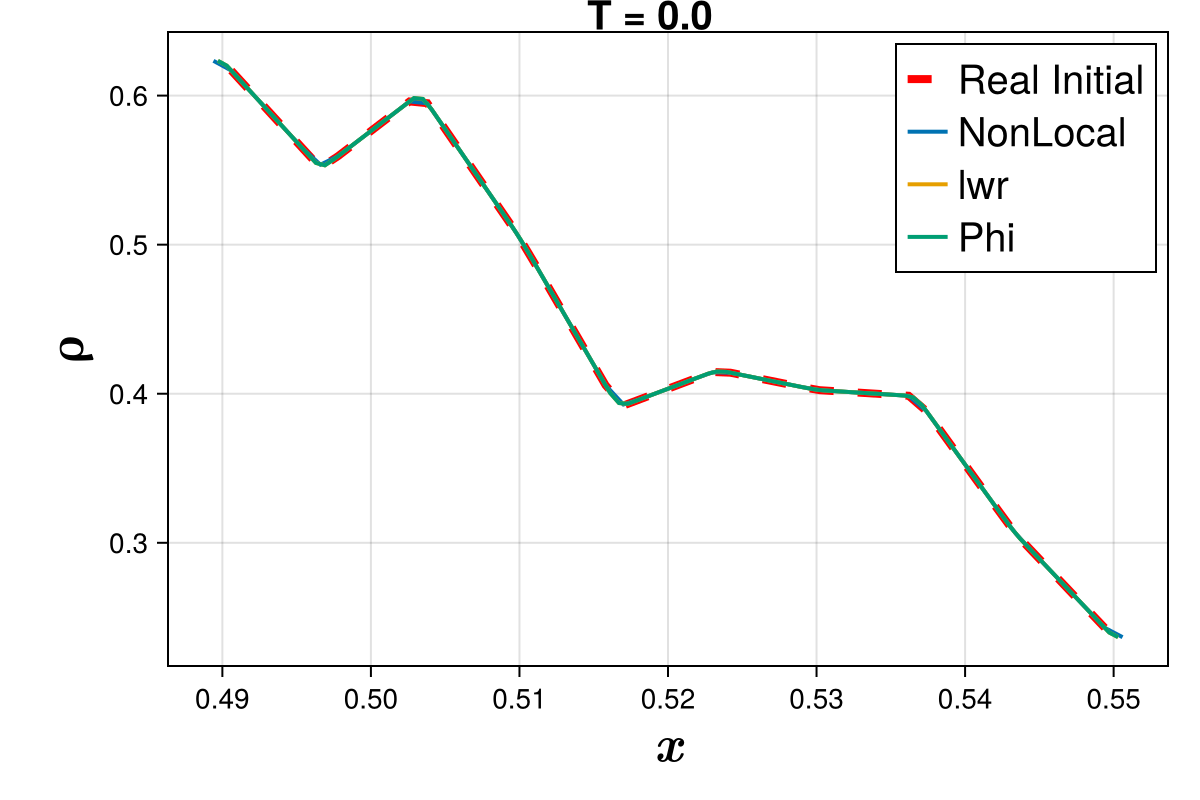}
    \includegraphics[width=0.32\linewidth]{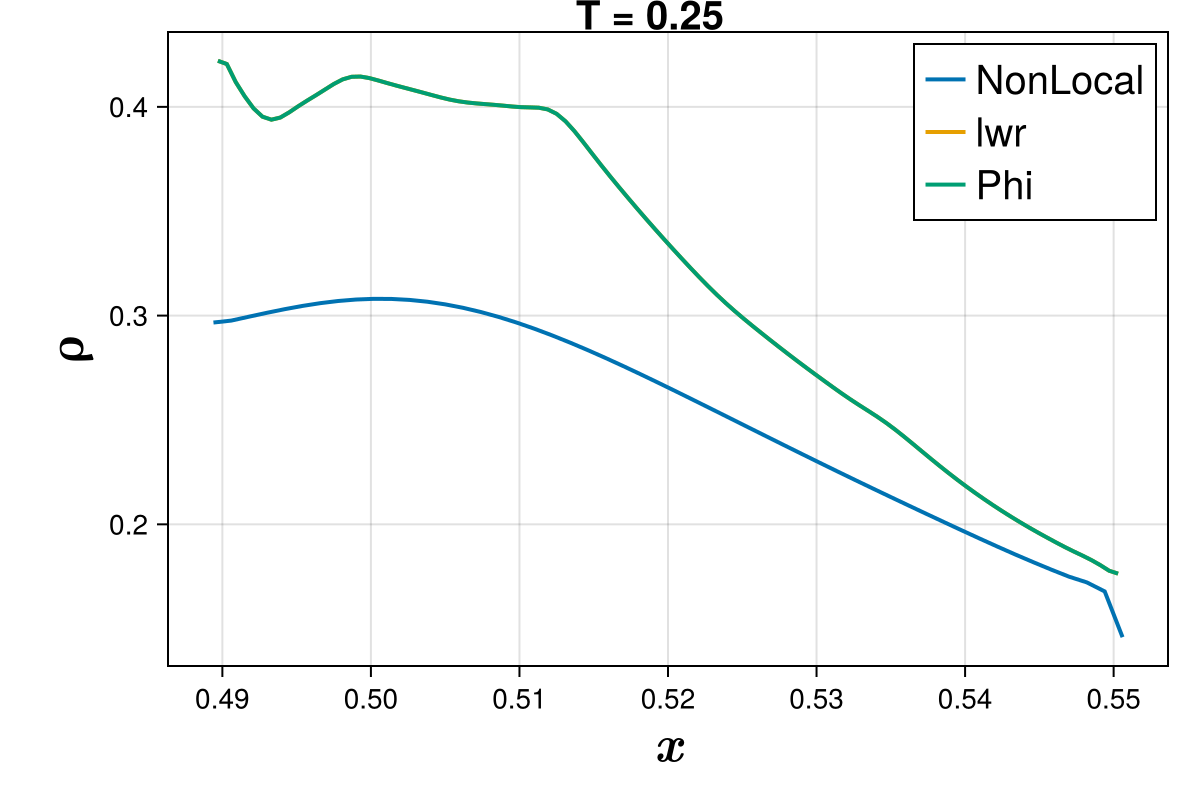}
    \includegraphics[width=0.32\linewidth]{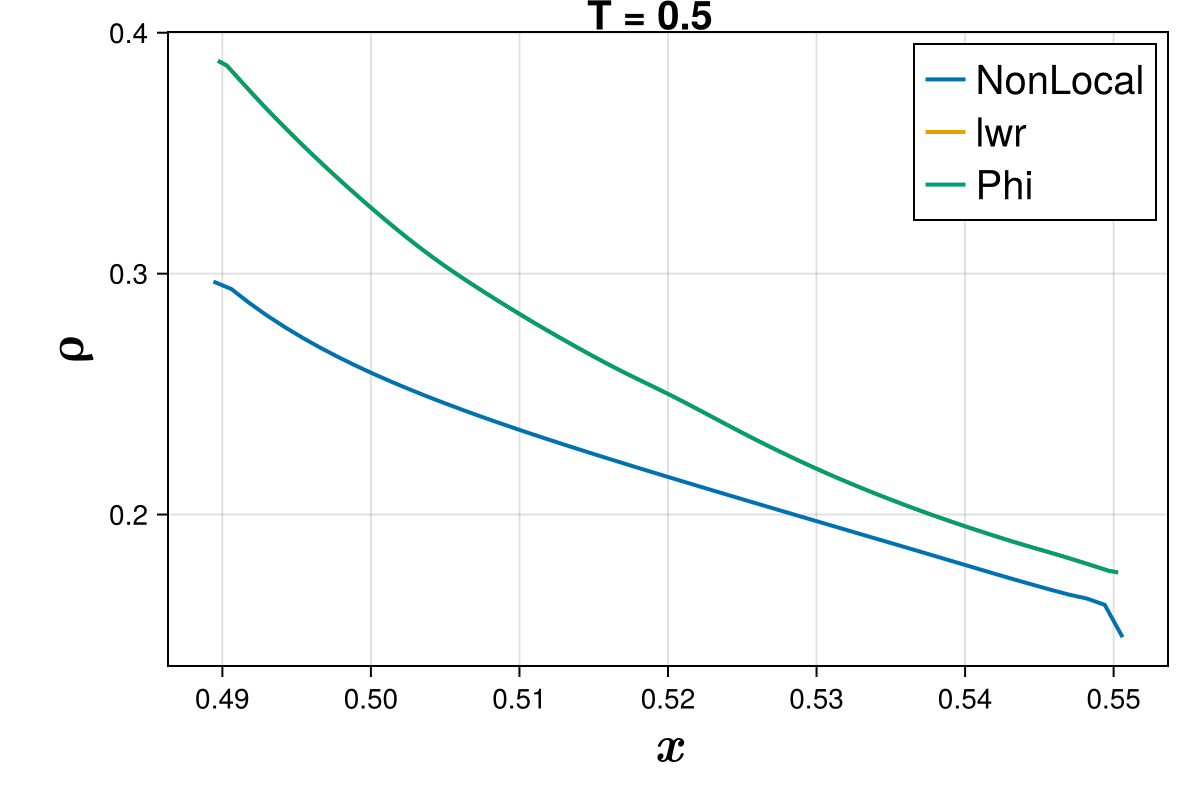}
    \includegraphics[width=0.32\linewidth]{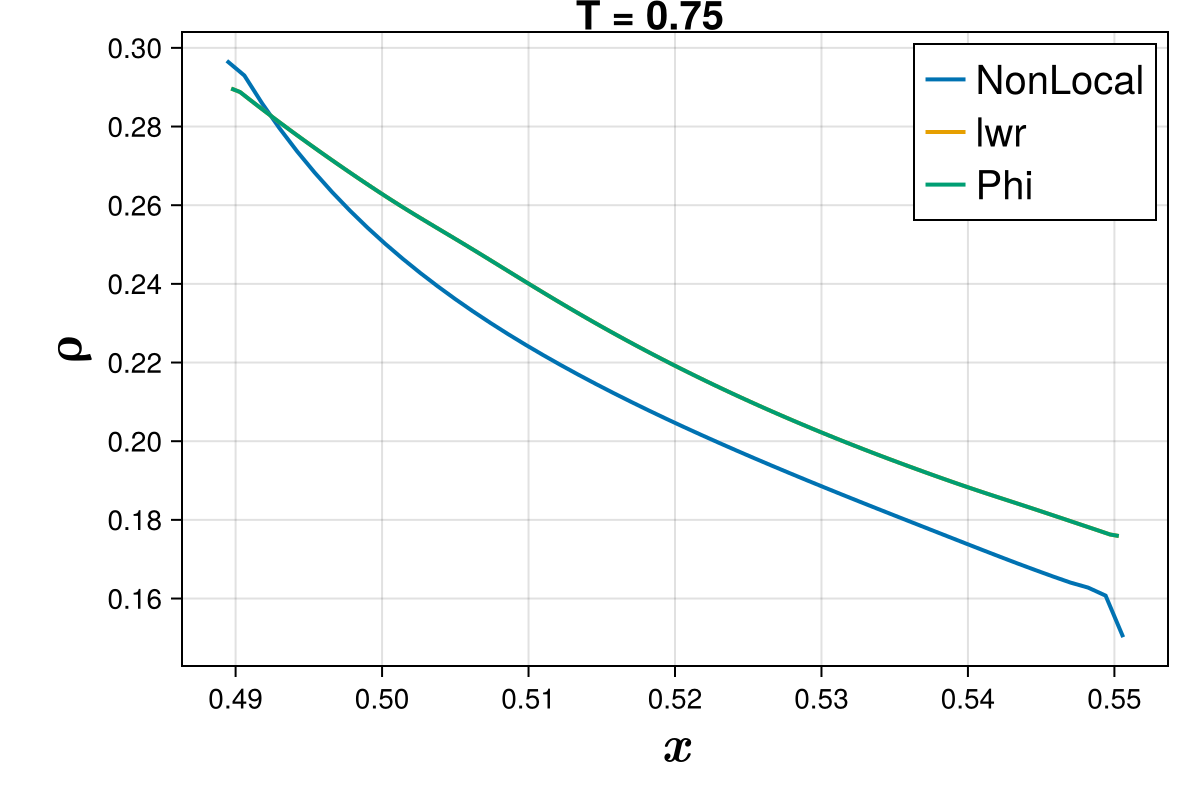}
    \includegraphics[width=0.32\linewidth]{Figures/Performance_final/dissipation_drone_gamma04_q4.png}
    \caption{\small{
    Step-by-step progression of performance in dissipation using the drone data. The frames correspond to simulation times of \( T = 0\),  \(0.25 \), \( 0.5 \), \( 0.75 \), and \( 1.0 \) hr, where \( T \) represents different time steps. The ground truth is the final density value from the data at \( T_G = 1.0 \) hr.}}

    \label{fig:dissipation_drone_evulotion}
\end{figure}

\section{Conclusion}
In this work we introduced a novel nonlocal traffic flow model to describe the macroscopic behavior of traffic flow dynamics. We introduced and motivated two alternative methods for calibrating the flow-density function. The first method for calibrating the FD that conceptually aligns with macroscopic flow theory (i.e. flow as a random process with a given probability transition kernel) and show that this method is successful in recreating empirical flow samples. We introduce a second calibration method that evaluates the model's PDE solution, selecting parameters that minimize the error between the model and empirical measurements. With this calibration, we show that the proposed model outperforms existing models by up to 80\%. The results of both calibrations provide evidence that underlying nonlocal dynamics \emph{and} diffusive correction are present in driver behavior. In each experiment, the best-performing models have non-zero values for the respective parameters. 
\appendix
\section{\textbf{Analytical preliminaries}} \label{S:prelim}
We start with two simple definitions from theory of stochastic processes, probability transition which provides the probability of transition between the states of density, and transition semigroup which explains the transition between the states of random variable (here the density function) during a certain period of time. For more detailed treatment of stochastic processes, we refer to classical textbooks on the topic (see e.g. \cite{dynkin1965markov, khasminskii2012stochastic}). For basics on measure theoretic probability theory, we refer the readers to \cite{feller1991introduction}. 
\begin{definition}[Probability transition] \label{def:prob_kernel}
    Let $(A, \mathcal A)$ be a measurable space. A map $p: A \times \mathcal A \to [0, 1]$ is called a probability transition kernel if 
    \begin{itemize}
        \item For any $\xi \in A$, the map $\mathcal A \ni B \mapsto p(\xi, A)$ is a probability measure on $(A, \mathcal A)$ 
        \item For every $B \in \mathcal A$, the map $A \ni \xi \mapsto p(\xi, A)$ is measurable. 
    \end{itemize}
\end{definition}
Probability transition can be considered as the probability of starting from a point $\xi$ and hitting a measurable set $A$. In the next definition, we include the temporal evolution into the definition of the probability transition. 
\begin{definition}[Transition Semigroup] \label{def:semigroup}
A collection $(p(t, \cdot, \cdot))_{t \in \RR_+}$ of probability transition kernels on $(A, \mathcal A)$, defines in the Definition \ref{def:prob_kernel}, is called a transition semigroup, if 
\begin{itemize}
    \item For any $\xi \in A$ and any $s \in \RR_+$, $p(0, \xi, d\zeta) = \delta_{\xi}(d\zeta)$ (i.e. Dirac at $\xi$). This in particular explains that the probability of starting from a point $\xi$ at a time $s$ and moving to a different point at the same time is negligible. 
    \item \tup{Chapman-Kolmogorov identity}: For any $s, t \ge 0$, and $B \in \mathcal A$, 
        \begin{equation*}
            p(t + s,\xi, B) = \int_A  p(s, \xi, dy) p(t, B) 
        \end{equation*}
        This property implies the uniqueness in the path of $t \mapsto p(t, \xi, d \zeta)$ which is a key point in proving that the claimed result is unique. 
    \item For every $B \in \mathcal A$, the function $(t, \xi) \mapsto p(t, \xi, B)$ is measurable with respect to the $\sigma$-field $\Bor(\RR_+) \otimes \mathcal A$. 
\end{itemize}
\end{definition}
Transition semigroup $p(\tau,x, B)$ can be simply interpreted as the probability of the associated random processes starting from time $t = 0$ at the position $x$ and then hitting a set $B$ at time $t = \tau$, i.e. after $\tau$ units of time. In addition, for any bounded measurable function $f$ on $\RR$, the collection of bounded linear operators $(T_t)_{t \ge 0}$ defined by 
\begin{equation*}
    T_t f(\xi) \Def \int p(t,  \xi, d \zeta) f(\zeta), \quad \xi \in \RR
\end{equation*}
is equivalently referred to as transition semigroup. 
\section{Kolmogorov Extension Theorem}\label{S:Kolmo}
In this section, we briefly discuss a key result in proving the existence and uniqueness of the probability measure on the space of all canonical maps. 
\begin{definition}[Consistent probability measures] \label{def:consistent} Let $\mathcal E(\RR_+)$ be the collection of all finite subsets of $\RR_+$. The collection of measures $\set{\mu_U: U \in \mathcal E(\RR_+)}$ is called consistent, if for any $U, V \in \mathcal E(\RR_+)$ such that $U \subset V$, $\mu_V\mid_{U} = \mu_U$. 
\end{definition}
\begin{theorem}[Kolmogorov Extension, \cite{oksendal2013stochastic}] Let $\mu_U$ be probability measure on $E^U$ such that the collection $\set{\mu_U: U \in \mathcal E(\RR_+)}$ be consistent in the sense of Definition \ref{def:consistent}. Then there exists a unique probability measure $\mu$ on $(\mathbb F(E), \mathcal F^*)$ such that $\mu \mid_{U} = \mu_U$.   
\end{theorem}

\bibliographystyle{siam} 
\bibliography{references.bib}

\end{document}